\input amstex



\def\spaces{\space\space\space\space\space\space\space\space\space\space}
\def\spacess{\message{\spaces\spaces\spaces\spaces\spaces\spaces\spaces}}
\spacess
\spacess
\message{Annals of Mathematics Style: Current Version: 1.1. June 10, 1992}
\spacess
\spacess

\catcode`\@=11

\hyphenation{acad-e-my acad-e-mies af-ter-thought anom-aly anom-alies
an-ti-deriv-a-tive an-tin-o-my an-tin-o-mies apoth-e-o-ses apoth-e-o-sis
ap-pen-dix ar-che-typ-al as-sign-a-ble as-sist-ant-ship as-ymp-tot-ic
asyn-chro-nous at-trib-uted at-trib-ut-able bank-rupt bank-rupt-cy
bi-dif-fer-en-tial blue-print busier busiest cat-a-stroph-ic
cat-a-stroph-i-cally con-gress cross-hatched data-base de-fin-i-tive
de-riv-a-tive dis-trib-ute dri-ver dri-vers eco-nom-ics econ-o-mist
elit-ist equi-vari-ant ex-quis-ite ex-tra-or-di-nary flow-chart
for-mi-da-ble forth-right friv-o-lous ge-o-des-ic ge-o-det-ic geo-met-ric
griev-ance griev-ous griev-ous-ly hexa-dec-i-mal ho-lo-no-my ho-mo-thetic
ideals idio-syn-crasy in-fin-ite-ly in-fin-i-tes-i-mal ir-rev-o-ca-ble
key-stroke lam-en-ta-ble light-weight mal-a-prop-ism man-u-script
mar-gin-al meta-bol-ic me-tab-o-lism meta-lan-guage me-trop-o-lis
met-ro-pol-i-tan mi-nut-est mol-e-cule mono-chrome mono-pole mo-nop-oly
mono-spline mo-not-o-nous mul-ti-fac-eted mul-ti-plic-able non-euclid-ean
non-iso-mor-phic non-smooth par-a-digm par-a-bol-ic pa-rab-o-loid
pa-ram-e-trize para-mount pen-ta-gon phe-nom-e-non post-script pre-am-ble
pro-ce-dur-al pro-hib-i-tive pro-hib-i-tive-ly pseu-do-dif-fer-en-tial
pseu-do-fi-nite pseu-do-nym qua-drat-ics quad-ra-ture qua-si-smooth
qua-si-sta-tion-ary qua-si-tri-an-gu-lar quin-tes-sence quin-tes-sen-tial
re-arrange-ment rec-tan-gle ret-ri-bu-tion retro-fit retro-fit-ted
right-eous right-eous-ness ro-bot ro-bot-ics sched-ul-ing se-mes-ter
semi-def-i-nite semi-ho-mo-thet-ic set-up se-vere-ly side-step sov-er-eign
spe-cious spher-oid spher-oid-al star-tling star-tling-ly
sta-tis-tics sto-chas-tic straight-est strange-ness strat-a-gem strong-hold
sum-ma-ble symp-to-matic syn-chro-nous topo-graph-i-cal tra-vers-a-ble
tra-ver-sal tra-ver-sals treach-ery turn-around un-at-tached un-err-ing-ly
white-space wide-spread wing-spread wretch-ed wretch-ed-ly Brown-ian
Eng-lish Euler-ian Feb-ru-ary Gauss-ian Grothen-dieck Hamil-ton-ian
Her-mit-ian Jan-u-ary Japan-ese Kor-te-weg Le-gendre Lip-schitz
Lip-schitz-ian Mar-kov-ian Noe-ther-ian No-vem-ber Rie-mann-ian
Schwarz-schild Sep-tem-ber Za-mo-lod-chi-kov Knizh-nik quan-tum Op-dam
Mac-do-nald Ca-lo-ge-ro Su-ther-land Mo-ser Ol-sha-net-sky  Pe-re-lo-mov
in-de-pen-dent ope-ra-tors
}

\Invalid@\nofrills
\Invalid@\usualspace
\newif\ifnofrills@
\def\nofrills@#1#2{\relaxnext@
  \DN@{\ifx\next\nofrills
    \nofrills@true\let#2\relax\DN@\nofrills{\nextii@}%
  \else
    \nofrills@false\def#2{#1}\let\next@\nextii@\fi
\next@}}
\def\usualspace@#1{\ifnofrills@\def\usualspace{#1}\fi}
\def\addto#1#2{\csname \expandafter\eat@\string#1@\endcsname
  \expandafter{\the\csname \expandafter\eat@\string#1@\endcsname#2}}
\newdimen\bigsize@
\def\big@#1#2{{\hbox{$\left#2\vcenter to#1\bigsize@{}%
  \right.\nulldelimiterspace\z@\m@th$}}}
\def\big{\big@\@ne}
\def\Big{\big@{1.5}}
\def\bigg{\big@\tw@}
\def\Bigg{\big@{2.5}}
\def\raggedcenter@{\leftskip\z@ plus.4\hsize \rightskip\leftskip
 \parfillskip\z@ \parindent\z@ \spaceskip.3333em \xspaceskip.5em
 \pretolerance9999\tolerance9999 \exhyphenpenalty\@M
 \hyphenpenalty\@M \let\\\linebreak}
\def\upperspecialchars{\def\ss{SS}\let\i=I\let\j=J\let\ae\AE\let\oe\OE
  \let\o\O\let\aa\AA\let\l\L}
\def\uppercasetext@#1{%
  {\spaceskip1.2\fontdimen2\the\font plus1.2\fontdimen3\the\font
   \upperspecialchars\uctext@#1$\m@th\aftergroup\eat@$}}
\def\uctext@#1$#2${\endash@#1-\endash@$#2$\uctext@}
\def\endash@#1-#2\endash@{%
\uppercase{#1}\if\notempty{#2}--\endash@#2\endash@\fi}
\def\runaway@#1{\DN@{#1}\ifx\envir@\next@
  \Err@{You seem to have a missing or misspelled \string\end#1 ...}%
  \let\envir@\empty\fi}
\newif\iftemp@
\def\notempty#1{TT\fi\def\test@{#1}\ifx\test@\empty\temp@false
  \else\temp@true\fi \iftemp@}

\font@\tensmc=cmcsc10
\font@\sevenex=cmex7 
\font@\sevenit=cmti7
\font@\eightrm=cmr8 
\font@\sixrm=cmr6 
\font@\eighti=cmmi8     \skewchar\eighti='177 
\font@\sixi=cmmi6       \skewchar\sixi='177   
\font@\eightsy=cmsy8    \skewchar\eightsy='60 
\font@\sixsy=cmsy6      \skewchar\sixsy='60   
\font@\eightex=cmex8 %
\font@\eightbf=cmbx8 
\font@\sixbf=cmbx6   
\font@\eightit=cmti8 
\font@\eightsl=cmsl8 
\font@\eightsmc=cmcsc10 
\font@\eighttt=cmtt8 

\loadmsam
\loadmsbm
\loadeufm
\UseAMSsymbols

\def\penaltyandskip@#1#2{\relax\ifdim\lastskip<#2\relax\removelastskip
      \ifnum#1=\z@\else\penalty@#1\relax\fi\vskip#2%
  \else\ifnum#1=\z@\else\penalty@#1\relax\fi\fi}
\def\nobreak{\penalty\@M
  \ifvmode\def\penalty@{\let\penalty@\penalty\count@@@}%
  \everypar{\let\penalty@\penalty\everypar{}}\fi}
\let\penalty@\penalty

\def\block{\RIfMIfI@\nondmatherr@\block\fi
       \else\ifvmode\vskip\abovedisplayskip\noindent\fi
        $$\def\endblock{\par\egroup$$}\fi
  \vbox\bgroup\advance\hsize-2\indenti\noindent}
\def\endblock{\par\egroup}

\def\logo@{\baselineskip2pc \hbox to\hsize{\hfil\eightpoint Typeset by
 \AmSTeX}}




\font\elevensc=cmcsc10 scaled\magstephalf
\font\tensc=cmcsc10

\font\eightsc=cmcsc10 scaled800

\font\elevenrm=cmr10 scaled \magstephalf
\font\ninerm=cmr9
\font\eightrm=cmr8
\font\sixrm=cmr6
\font\fiverm=cmr5

\font\eleveni=cmmi10 scaled\magstephalf
\font\ninei=cmmi9
\font\eighti=cmmi8
\font\sixi=cmmi6
\font\fivei=cmmi5
\skewchar\ninei='177 \skewchar\eighti='177 \skewchar\sixi='177
\skewchar\eleveni='177

\font\elevensy=cmsy10 scaled\magstephalf
\font\ninesy=cmsy9
\font\eightsy=cmsy8
\font\sixsy=cmsy6
\font\fivesy=cmsy5
\skewchar\ninesy='60 \skewchar\eightsy='60 \skewchar\sixsy='60
\skewchar\elevensy'60

\font\eighteenbf=cmbx10 scaled\magstep3

\font\twelvebf=cmbx10 scaled \magstep1
\font\elevenbf=cmbx10 scaled \magstephalf
\font\tenbf=cmbx10
\font\ninebf=cmbx9
\font\eightbf=cmbx8
\font\sixbf=cmbx6
\font\fivebf=cmbx5

\font\elevenit=cmti10 scaled\magstephalf
\font\nineit=cmti9
\font\eightit=cmti8

\font\eighteenmib=cmmib10 scaled \magstep3
\font\twelvemib=cmmib10 scaled \magstep1
\font\elevenmib=cmmib10 scaled\magstephalf
\font\tenmib=cmmib10
\font\eightmib=cmmib10 scaled 800 
\font\sixmib=cmmib10 scaled 600

\font\eighteensyb=cmbsy10 scaled \magstep3 
\font\twelvesyb=cmbsy10 scaled \magstep1
\font\elevensyb=cmbsy10 scaled \magstephalf
\font\tensyb=cmbsy10 
\font\eightsyb=cmbsy10 scaled 800
\font\sixsyb=cmbsy10 scaled 600
 
\font\elevenex=cmex10 scaled \magstephalf
\font\tenex=cmex10     
\font\eighteenex=cmex10 scaled \magstep3


\def\elevenpoint{\def\rm{\fam0\elevenrm}%
  \textfont0=\elevenrm \scriptfont0=\eightrm \scriptscriptfont0=\sixrm
  \textfont1=\eleveni \scriptfont1=\eighti \scriptscriptfont1=\sixi
  \textfont2=\elevensy \scriptfont2=\eightsy \scriptscriptfont2=\sixsy
  \textfont3=\elevenex \scriptfont3=\tenex \scriptscriptfont3=\tenex
  \def\bf{\fam\bffam\elevenbf}%
  \def\it{\fam\itfam\elevenit}%
  \textfont\bffam=\elevenbf \scriptfont\bffam=\eightbf
   \scriptscriptfont\bffam=\sixbf
\normalbaselineskip=13.95pt
  \setbox\strutbox=\hbox{\vrule height9.5pt depth4.4pt width0pt\relax}%
  \normalbaselines\rm}

\elevenpoint 

\def\ninepoint{\def\rm{\fam0\ninerm}%
  \textfont0=\ninerm \scriptfont0=\sixrm \scriptscriptfont0=\fiverm
  \textfont1=\ninei \scriptfont1=\sixi \scriptscriptfont1=\fivei
  \textfont2=\ninesy \scriptfont2=\sixsy \scriptscriptfont2=\fivesy
  \textfont3=\tenex \scriptfont3=\tenex \scriptscriptfont3=\tenex
  \def\it{\fam\itfam\nineit}%
  \textfont\itfam=\nineit
  \def\bf{\fam\bffam\ninebf}%
  \textfont\bffam=\ninebf \scriptfont\bffam=\sixbf
   \scriptscriptfont\bffam=\fivebf
\normalbaselineskip=11pt
  \setbox\strutbox=\hbox{\vrule height8pt depth3pt width0pt\relax}%
  \normalbaselines\rm}

\def\eightpoint{\def\rm{\fam0\eightrm}%
  \textfont0=\eightrm \scriptfont0=\sixrm \scriptscriptfont0=\fiverm
  \textfont1=\eighti \scriptfont1=\sixi \scriptscriptfont1=\fivei
  \textfont2=\eightsy \scriptfont2=\sixsy \scriptscriptfont2=\fivesy
  \textfont3=\tenex \scriptfont3=\tenex \scriptscriptfont3=\tenex
  \def\it{\fam\itfam\eightit}%
  \textfont\itfam=\eightit
  \def\bf{\fam\bffam\eightbf}%
  \textfont\bffam=\eightbf \scriptfont\bffam=\sixbf
   \scriptscriptfont\bffam=\fivebf
\normalbaselineskip=12pt
  \setbox\strutbox=\hbox{\vrule height8.5pt depth3.5pt width0pt\relax}%
  \normalbaselines\rm}


\def\eighteenbold{\def\rm{\fam0\eighteenbf}%
  \textfont0=\eighteenbf \scriptfont0=\twelvebf \scriptscriptfont0=\tenbf
  \textfont1=\eighteenmib \scriptfont1=\twelvemib\scriptscriptfont1=\tenmib
  \textfont2=\eighteensyb \scriptfont2=\twelvesyb\scriptscriptfont2=\tensyb
  \textfont3=\eighteenex \scriptfont3=\tenex \scriptscriptfont3=\tenex
  \def\bf{\fam\bffam\eighteenbf}%
  \textfont\bffam=\eighteenbf \scriptfont\bffam=\twelvebf
   \scriptscriptfont\bffam=\tenbf
\normalbaselineskip=20pt
  \setbox\strutbox=\hbox{\vrule height13.5pt depth6.5pt width0pt\relax}%
\everymath {\fam0 }
\everydisplay {\fam0 }
  \normalbaselines\rm}

\def\elevenbold{\def\rm{\fam0\elevenbf}%
  \textfont0=\elevenbf \scriptfont0=\eightbf \scriptscriptfont0=\sixbf
  \textfont1=\elevenmib \scriptfont1=\eightmib \scriptscriptfont1=\sixmib
  \textfont2=\elevensyb \scriptfont2=\eightsyb \scriptscriptfont2=\sixsyb
  \textfont3=\elevenex \scriptfont3=\elevenex \scriptscriptfont3=\elevenex
  \def\bf{\fam\bffam\elevenbf}%
  \textfont\bffam=\elevenbf \scriptfont\bffam=\eightbf
   \scriptscriptfont\bffam=\sixbf
\normalbaselineskip=14pt
  \setbox\strutbox=\hbox{\vrule height10pt depth4pt width0pt\relax}%
\everymath {\fam0 }
\everydisplay {\fam0 }
  \normalbaselines\bf}

\hsize=31pc
\vsize=48pc

\parindent=22pt
\parskip=0pt

\widowpenalty=10000
\clubpenalty=10000

\topskip=12pt 

\skip\footins=20pt
\dimen\footins=3in 

\abovedisplayskip=6.95pt plus3.5pt minus 3pt
\belowdisplayskip=\abovedisplayskip


\voffset=7pt\hoffset= .7in

\newif\iftitle

\def\amheadline{\iftitle%
\hbox to\hsize{\hss\currannalsline\hss}\else\line{\ifodd\pageno
\hfill\thetitle\hfill\llap{\elevenrm\folio}\else\rlap{\elevenrm\folio}
\hfill\theauthors\hfill\fi}\fi}

\headline={\amheadline}
\footline={\global\titlefalse}


\def\annalsline#1#2{\vfill\eject
\ifodd\pageno\else 
\line{\hfill}
\vfill\eject\fi
\global\titletrue
\def\currannalsline{\eightrm 
{\eightbf#1} (#2), \thepages}}

\def\titleheadline#1{\def\one{#1}\ifx\one\empty\else
\def\thetitle{{
\let\\ \relax\eightsc\uppercase{#1}}}\fi}

\newif\ifshort

\let\shorttitle\titleheadline

\def\onpages#1#2{\def\thepages{#1--#2}}

\def\thismuchskip[#1]{\vskip#1pt}
\def\ilook{\ifx\next[ \let\go\thismuchskip\else
\let\go\relax\vskip1pt\fi\go}

\def\institution#1{\def\theinstitutions{\vbox{\baselineskip10pt
\def\and{{\eightrm and }}
\def\\{\futurelet\next\ilook}\eightsc #1}}}
\let\institutions\institution

\newwrite\auxfile

\def\startingpage#1{\def\one{#1}\ifx\one\empty\global\pageno=1\else
\global\pageno=#1\fi
\theoremcount=0 \eqcount=0 \sectioncount=0 
\openin1 \jobname.aux \ifeof1 
\onpages{#1}{???}
\else\closein1 \relax\input \jobname.aux
\onpages{#1}{\lastpage}
\fi\immediate\openout\auxfile=\jobname.aux
}

\def\endarticle{\ifRefsUsed\global\RefsUsedfalse%
\else\vskip21pt\theinstitutions%
\nobreak\vskip8pt
\fi%
\write\auxfile{\string\def\string\lastpage{\the\pageno}}}

\outer\def\bye{\endarticle\par \vfill \supereject \end}

\def\document{\let\fontlist@\relax\let\alloclist@\relax
 \elevenpoint}


\newif\ifacks
\long\def\acknowledgements#1{\def\one{#1}\ifx\one\empty\else
\vskip-\baselineskip%
\global\ackstrue\footnote{\ \unskip}{*#1}\fi}

\def\title#1{\titleheadline{#1}
\vbox to80pt{\vfill
\baselineskip=18pt
\parindent=0pt
\overfullrule=0pt
\hyphenpenalty=10000
\everypar={\hskip\parfillskip\relax}
\hbadness=10000
\def\\ {\vskip1sp}
\eighteenbold#1\vskip1sp}}

\newif\ifauthor

\def\author#1{\vskip11pt
\hbox to\hsize{\hss\tenrm By \tensc#1\ifacks\global\acksfalse*\fi\hss}
\ifshort\else\xdef\theauthors{{\eightsc\uppercase{#1}}}\fi%
\vskip21pt\global\authortrue\everypar={\global\authorfalse\everypar={}}}

\def\twoauthors#1#2{\vskip11pt
\hbox to\hsize{\hss%
\tenrm By \tensc#1 {\tenrm and} #2\ifacks\global\acksfalse*\fi\hss}
\ifshort\else\xdef\theauthors{{\eightsc\uppercase{#1 and #2}}}\fi%
\vskip21pt
\global\authortrue\everypar={\global\authorfalse\everypar={}}}


\newcount\theoremcount
\newcount\sectioncount
\newcount\eqcount

\newif\ifspecialnumon

\def\eqnumber=#1 {\global\eqcount=#1 \global\advance\eqcount by-1\relax}
\def\sectionnumber=#1 {\global\sectioncount=#1 
\global\advance\sectioncount by-1\relax}
\def\proclaimnumber=#1 {\global\theoremcount=#1 
\global\advance\theoremcount by-1\relax}

\newif\ifsection
\newif\ifsubsection

\def\elevenboldmath#1{$#1$\egroup}
\def\mathbold{\hbox\bgroup\elevenbold\elevenboldmath}

\def\section#1{\global\theoremcount=0
\global\eqcount=0
\ifauthor\global\authorfalse\else%
\vskip18pt plus 18pt minus 6pt\fi%
{\parindent=0pt
\everypar={\hskip\parfillskip}
\def\\ {\vskip1sp}\elevenpoint\bf%
\ifspecialnumon\global\specialnumonfalse$\rm\spnum$%
\gdef\sectnum{$\rm\spnum$}%
\else\interlinepenalty=10000%
\global\advance\sectioncount by1\relax\the\sectioncount%
\gdef\sectnum{\the\sectioncount}%
\fi. \hskip6pt#1
\vrule width0pt depth12pt}
\hskip\parfillskip
\global\sectiontrue%
\everypar={\global\sectionfalse\global\interlinepenalty=0\everypar={}}%
\ignorespaces

}


\newif\ifspequation

\let\eqno\leqno 

\newif\ifineqalignno
\let\saveleqalignno\leqalignno                        
\def\leqalignno{\let\eqnu\Eeqnu\saveleqalignno}

\let\eqalignno\leqalignno

\def\sectandeqnum{%
\ifspecialnumon\global\specialnumonfalse
$\rm\spnum$\gdef\eqnum{{$\rm\spnum$}}\else\global\firstlettertrue
\global\advance\eqcount by1 
\ifappend\applett\else\the\sectioncount\fi.%
\the\eqcount
\xdef\eqnum{\ifappend\applett\else\the\sectioncount\fi.\the\eqcount}\fi}

\def\eqnu{\leqno{\hbox{\elevenrm\ifspequation\else(\fi\sectandeqnum
\ifspequation\global\spequationfalse\else)\fi}}}      

\def\Speqnu{\global\setbox\leqnobox=\hbox{\elevenrm
\ifspequation\else%
(\fi\sectandeqnum\ifspequation\global\spequationfalse\else)\fi}}

\def\Eeqnu{\hbox{\elevenrm
\ifspequation\else%
(\fi\sectandeqnum\ifspequation\global\spequationfalse\else)\fi}}

\newif\iffirstletter
\global\firstlettertrue
\def\eqletter#1{\global\specialnumontrue\iffirstletter\global\firstletterfalse
\global\advance\eqcount by1\fi
\gdef\spnum{\the\sectioncount.\the\eqcount#1}\eqnu}

\newbox\leqnobox
\def\outsideeqnu#1{\global\setbox\leqnobox=\hbox{#1}}

\def\eatone#1{}

\def\dosplit#1#2{\vskip-.5\abovedisplayskip
\setbox0=\hbox{$\let\eqno\outsideeqnu%
\let\eqnu\Speqnu\let\leqno\outsideeqnu#2$}%
\setbox1\vbox{\noindent\hskip\wd\leqnobox\ifdim\wd\leqnobox>0pt\hskip1em\fi%
$\displaystyle#1\mathstrut$\hskip0pt plus1fill\relax
\vskip1pt
\line{\hfill$\let\eqnu\eatone\let\leqno\eatone%
\displaystyle#2\mathstrut$\ifmathqed~~\qed\fi}}%
\copy1
\ifvoid\leqnobox
\else\dimen0=\ht1 \advance\dimen0 by\dp1
\vskip-\dimen0
\vbox to\dimen0{\vfill
\hbox{\unhbox\leqnobox}
\vfill}
\fi}

\everydisplay{\lookforbreak}

\long\def\lookforbreak #1$${\def\mathone{#1}
\expandafter\testforbreak\mathone\splitmath @}

\def\testforbreak#1\splitmath #2@{\def\mathtwo{#2}\ifx\mathtwo\empty%
#1$$%
\ifmathqed\vskip-\belowdisplayskip
\setbox0=\vbox{\let\eqno\relax\let\eqnu\relax$\displaystyle#1$}%
\vskip-\ht0\vskip-3.5pt\hbox to\hsize{\hfill\qed}
\vskip\ht0\vskip3.5pt\fi
\else$$\vskip-\belowdisplayskip
\vbox{\dosplit{#1}{\let\eqno\eatone
\let\splitmath\relax#2}}%
\nobreak\vskip.5\belowdisplayskip
\noindent\ignorespaces\fi}


\newif\ifmathqed



\newcount\linenum
\newcount\colnum

\def\spline{\omit&\multispan{\the\colnum}{\hrulefill}\cr}
\def\colcounter{\ifnum\linenum=1\global\advance\colnum by1\fi}

\def\everyline{\noalign{\global\advance\linenum by1\relax}%
\ifnum\linenum=2\spline\fi}

\def\mtable{\bgroup\offinterlineskip
\everycr={\everyline}\global\linenum=0
\halign\bgroup\vrule height 10pt depth 4pt width0pt
\hfill$##$\hfill\hskip6pt\ifnum\linenum>1
\vrule\fi&&\colcounter\hskip12pt\hfill$##$\hfill\hskip12pt\cr}

\def\endmtable{\crcr\egroup\egroup}




\def\xast{*}
\newcount\intable
\newcount\mathcol
\newcount\savemathcol
\newcount\topmathcol
\newdimen\arrayhspace
\newdimen\arrayvspace

\arrayhspace=8pt 
\arrayvspace=12pt 

\newif\ifdollaron

\def\mathalign#1{\def\arg{#1}\ifx\arg\xast%
\let\go\relax\else\let\go\mathalign%
\global\advance\mathcol by1 %
\global\advance\topmathcol by1 %
\expandafter\def\csname  mathcol\the\mathcol\endcsname{#1}%
\fi\go}

\def\arraypickapart#1]#2*{\if#1c \ifmmode\vcenter\else
\global\dollarontrue$\vcenter\fi\else%
\if#1t\vtop\else\if#1b\vbox\fi\fi\fi\bgroup%
\def\one{#2}}

\def\arraystrut{\vrule height .7\arrayvspace depth .3\arrayvspace width 0pt}

\def\array#1#2*{\def\firstarg{#1}%
\if\firstarg[ \def\two{#2} \expandafter\arraypickapart\two*\else%
\ifmmode\vcenter\else\vbox\fi\bgroup \def\one{#1#2}\fi%
\global\everycr={\noalign{\global\mathcol=\savemathcol\relax}}%
\def\\ {\cr}%
\global\advance\intable by1 %
\ifnum\intable=1 \global\mathcol=0 \savemathcol=0 %
\else \global\advance\mathcol by1 \savemathcol=\mathcol\fi%
\expandafter\mathalign\one*%
\mathcol=\savemathcol %
\halign\bgroup&\hskip.5\arrayhspace\arraystrut%
\global\advance\mathcol by1 \relax%
\expandafter\if\csname mathcol\the\mathcol\endcsname r\hfill\else%
\expandafter\if\csname mathcol\the\mathcol\endcsname c\hfill\fi\fi%
$\displaystyle##$%
\expandafter\if\csname mathcol\the\mathcol\endcsname r\else\hfill\fi\relax%
\hskip.5\arrayhspace\cr}

\def\endarray{\crcr\egroup\egroup%
\global\mathcol=\savemathcol %
\global\advance\intable by -1\relax%
\ifnum\intable=0 %
\ifdollaron\global\dollaronfalse $\fi
\loop\ifnum\topmathcol>0 %
\expandafter\def\csname  mathcol\the\topmathcol\endcsname{}%
\global\advance\topmathcol by-1 \repeat%
\global\everycr={}\fi%
}

\def\big#1{{\hbox{$\left#1\vbox to 10pt{}\right.\n@space$}}}
\def\Big#1{{\hbox{$\left#1\vbox to 13pt{}\right.\n@space$}}}
\def\bigg#1{{\hbox{$\left#1\vbox to 16pt{}\right.\n@space$}}}
\def\Bigg#1{{\hbox{$\left#1\vbox to 19pt{}\right.\n@space$}}}


\def\figcaption#1#2#3{\topinsert
\vskip4pt 
\vbox to#3{\vfill}\vskip1sp
\setbox0=\hbox{\eightsc Figure #1.\hskip12pt\eightpoint #2}
\ifdim\wd0>\hsize
\noindent\eightsc Figure #1.\hskip12pt\eightpoint #2
\else
\centerline{\eightsc Figure #1.\hskip12pt\eightpoint #2}
\fi
\vskip16pt
\endinsert}

\def\wfig#1#2#3{\topinsert
\vskip4pt 
\hbox to\hsize{\hss\vbox{\hrule height .25pt width #3
\hbox to #3{\vrule width .25pt height #2\hfill\vrule width .25pt height #2}
\hrule height.25pt}\hss}
\vskip1sp
\centerline{\eightsc Figure #1}
\vskip16pt
\endinsert}

\def\wfigcaption#1#2#3#4{\topinsert
\vskip4pt 
\hbox to\hsize{\hss\vbox{\hrule height .25pt width #4
\hbox to #4{\vrule width .25pt height #3\hfill\vrule width .25pt height #3}
\hrule height.25pt}\hss}
\vskip1sp
\setbox0=\hbox{\eightsc Figure #1.\hskip12pt\eightpoint\rm #2}
\ifdim\wd0>\hsize
\noindent\eightsc Figure #1.\hskip12pt\eightpoint\rm #2\else
\centerline{\eightsc Figure #1.\hskip12pt\eightpoint\rm #2}\fi
\vskip16pt
\endinsert}

\def\tabcaption#1#2{\vskip6pt
\setbox0=\hbox{\eightsc Table #1.\hskip12pt\eightpoint #2}
\ifdim\wd0>\hsize
\noindent\eightsc Table #1.\hskip12pt\eightpoint #2
\else
\centerline{\eightsc Table #1.\hskip12pt\eightpoint #2}
\fi
\vskip6pt}

\def\endinsert{\egroup\if@mid\dimen@\ht\z@\advance\dimen@\dp\z@ 
\advance\dimen@ 12\p@\advance\dimen@\pagetotal\ifdim\dimen@ >\pagegoal 
\@midfalse\p@gefalse\fi\fi\if@mid\smallskip\box\z@\bigbreak\else
\insert\topins{\penalty 100 \splittopskip\z@skip\splitmaxdepth\maxdimen
\floatingpenalty\z@\ifp@ge\dimen@\dp\z@\vbox to\vsize {\unvbox \z@ 
\kern -\dimen@ }\else\box\z@\nobreak\smallskip\fi}\fi\endgroup}

\def\pagecontents{
\ifvoid\topins \else\iftitle\else 
\unvbox \topins \fi\fi \dimen@ =\dp \@cclv \unvbox 
\@cclv 
\ifvoid\topins\else\iftitle\unvbox\topins\fi\fi
\ifvoid \footins \else \vskip \skip \footins \footnoterule 
\unvbox \footins \fi \ifr@ggedbottom \kern -\dimen@ \vfil \fi}


\newif\ifappend

\def\appendix#1#2{\def\applett{#1}\def\two{#2}%
\global\appendtrue
\global\theoremcount=0
\global\eqcount=0
\vskip18pt plus 18pt
\vbox{\parindent=0pt
\everypar={\hskip\parfillskip}
\def\\ {\vskip1sp}\elevenbold Appendix%
\ifx\applett\empty\gdef\applett{A}\ifx\two\empty\else.\fi%
\else\ #1.\fi\hskip6pt#2\vskip12pt}%
\global\sectiontrue%
\everypar={\global\sectionfalse\everypar={}}\nobreak\ignorespaces}

\newif\ifRefsUsed
\long\def\references{\global\RefsUsedtrue\vskip21pt
\theinstitutions
\global\everypar={}\global\bibnum=0
\vskip20pt\goodbreak\bgroup
\vbox{\centerline{\eightsc References}\vskip6pt}%
\ifdim\maxbibwidth>0pt
\leftskip=\maxbibwidth%
\parindent=-\maxbibwidth%
\else
\leftskip=18pt%
\parindent=-18pt%
\fi
\ninepoint
\frenchspacing
\nobreak\ignorespaces\everypar={\amref}%
}

\def\endreferences{\vskip1sp\egroup\global\everypar={}%
\nobreak\vskip8pt\vbox{\thereceived\therevised}
}

\newcount\bibnum

\def\amref#1 {\global\advance\bibnum by1%
\immediate\write\auxfile{\string\expandafter\string\def\string\csname
\space #1croref\string\endcsname{[\the\bibnum]}}%
\leavevmode\hbox to18pt{\hbox to13.2pt{\hss[\the\bibnum]}\hfill}}

\def\bibline{\hbox to30pt{\hrulefill}\/\/}

\def\name#1{{\eightsc#1}}

\newdimen\maxbibwidth
\def\AuthorRefNames [#1] {%
\immediate\write\auxfile{\string\def\string\cite\string##1{[\string##1]}}

\def\amref{\spamref}
\setbox0=\hbox{[#1] }\global\maxbibwidth=\wd0\relax}

\def\spamref[#1] {\leavevmode\hbox to\maxbibwidth{\hss[#1]\hfill}}


\def\footnoterule{\kern-3pt\hrule width1in height.5pt\kern2.5pt}

\def\footnote#1#2{%
\plainfootnote{#1}{{\eightpoint\normalbaselineskip11pt
\normalbaselines#2}}}

\def\vfootnote#1{%
\insert \footins \bgroup \eightpoint\baselineskip11pt
\interlinepenalty \interfootnotelinepenalty
\splittopskip \ht \strutbox \splitmaxdepth \dp \strutbox \floatingpenalty 
\@MM \leftskip \z@skip \rightskip \z@skip \spaceskip \z@skip 
\xspaceskip \z@skip
{#1}$\,$\footstrut \futurelet \next \fo@t}


\newif\iffirstadded
\newif\ifadded

\def\addedlett{}

\def\alltheoremnums{%
\ifspecialnumon\global\specialnumonfalse
\ifadded\global\addedfalse
\iffirstadded\global\firstaddedfalse
\global\advance\theoremcount by1 \fi
\ifappend\applett\else\the\sectioncount\fi.\the\theoremcount\addedlett%
\xdef\theoremnum{\ifappend\applett\else\the\sectioncount\fi.%
\the\theoremcount\addedlett}%
\else$\rm\spnum$\def\theoremnum{{$\rm\spnum$}}\fi%
\else\global\firstaddedtrue
\global\advance\theoremcount by1 
\ifappend\applett\else\the\sectioncount\fi.\the\theoremcount%
\xdef\theoremnum{\ifappend\applett\else\the\sectioncount\fi.%
\the\theoremcount}\fi}

\def\allcorolnums{%
\ifspecialnumon\global\specialnumonfalse
\ifadded\global\addedfalse
\iffirstadded\global\firstaddedfalse
\global\advance\corolcount by1 \fi
\the\corolcount\addedlett%
\else$\rm\spnum$\def\corolnum{$\rm\spnum$}\fi%
\else\global\advance\corolcount by1 
\the\corolcount\fi}


\newcount\corolcount
\def\xcorol{Corollary}
\def\xtheorem{Theorem}
\def\xmaintheorem{Main Theorem}

\newif\ifthtitle

\let\saverparen)
\let\savelparen(
\def\rmparenl{{\rm(}}
\def\rmparenr{{\rm\/)}}
{
\catcode`(=13
\catcode`)=13
\gdef\makeparensRM{\catcode`(=13\catcode`)=13\let(=\rmparenl%
\let)=\rmparenr%
\everymath{\let(\savelparen%
\let)\saverparen}%
\everydisplay{\let(\savelparen%
\let)\saverparen\lookforbreak}}}

\medskipamount=8pt plus.1\baselineskip minus.05\baselineskip

\def\rmtext#1{\hbox{\rm#1}}

\def\proclaim#1{\vskip-\lastskip
\def\one{#1}\ifx\one\xtheorem\global\corolcount=0\fi
\ifsection\global\sectionfalse\vskip-6pt\fi
\medskip
{\elevensc#1}%
\ifx\one\xmaintheorem\global\corolcount=0
\gdef\theoremnum{Main Theorem}\else%
\ifx\one\xcorol\ 
\alltheoremnums 
\else\ \alltheoremnums\fi\fi%
\ifthtitle\ \global\thtitlefalse{\rm(\thethtitle)}\fi.%
\hskip1em\bgroup\let\text\rmtext\makeparensRM\it\ignorespaces}

\def\nonumproclaim#1{\vskip-\lastskip
\def\one{#1}\ifx\one\xtheorem\global\corolcount=0\fi
\ifsection\global\sectionfalse\vskip-6pt\fi
\medskip
{\elevensc#1}.\ifx\one\xmaintheorem\global\corolcount=0
\gdef\theoremnum{Main Theorem}\fi\hskip.5pc%
\bgroup\it\makeparensRM\ignorespaces}

\def\endproclaim{\egroup\medskip}


\def\xproof{Proof}
\def\xremark{Remark}
\def\xcase{Case}
\def\xsubcase{Subcase}
\def\xconjecture{Conjecture}
\def\xstep{Step}
\def\xof{of}

\def\deconstruct#1 #2 #3 #4 #5 @{\def\one{#1}\def\two{#2}\def\three{#3}%
\def\four{#4}%
\ifx\two\empty #1\else%
\ifx\one\xproof%
\ifx\two\xof%
  \ifx\three\xcorol Proof of Corollary \rm#4\else%
     \ifx\three\xtheorem Proof of Theorem \rm#4\else\xone\fi%
  \fi\fi%
\else\xone\fi\fi.}

\def\pickup#1 {\def\this{#1}%
\ifx\this\xproof\global\let\go\demoproof
\global\let\enddemo\endproof\else
\ifx\this\xremark\global\let\go\demoremark\else
\ifx\this\xcase\global\let\go\demostep\else
\ifx\this\xsubcase\global\let\go\demostep\else
\ifx\this\xconjecture\global\let\go\demostep\else
\ifx\this\xstep\global\let\go\demostep\else
\global\let\go\demoproof\fi\fi\fi\fi\fi\fi}

\newif\ifnonum
\def\demo#1{\vskip-\lastskip
\ifsection\global\sectionfalse\vskip-6pt\fi
\def\one{#1 }\def\two{#1*}%
\setbox0=\hbox{\expandafter\pickup\one}\expandafter\go\two}

\def\numbereddemo#1{\vskip-\lastskip
\ifsection\global\sectionfalse\vskip-6pt\fi
\def\two{#1*}%
\expandafter\demoremark\two}

\def\demoproof#1*{\medskip\def\xone{#1}
{\ignorespaces\it\expandafter\deconstruct\xone {} {} {} {} {} @%
\unskip\hskip6pt}\rm\ignorespaces}

\def\demoremark#1*{\medskip
{\it\ignorespaces#1\/} \ifnonum\global\nonumtrue\else
 \alltheoremnums\unskip.\fi\hskip1pc\rm\ignorespaces}

\def\demostep#1 #2*{\vskip4pt
{\it\ignorespaces#1\/} #2.\hskip1pc\rm\ignorespaces}

\def\enddemo{\medskip}

\def\endproof{\ifmathqed\global\mathqedfalse\medskip\else
\parfillskip=0pt~~\hfill\qed\medskip
\fi\global\parfillskip0pt plus 1fil\relax
\gdef\enddemo{\medskip}}

\def\qed{\vbox{\hrule\hbox{\vrule height6pt\hskip6pt\vrule}\hrule}}








\def\stripbs#1#2*{\def\one{#2}}

\def\emptyspace{ }
\def\nextthing{}
\def\newline{***}
\def\eatone#1{ }

\def\lookatspace#1{\ifcat\noexpand#1\ \else%
\gdef\nextthing{}\xdef\next{#1}%
\ifx\next\emptyspace%
\let\nextthing\emptyspace\else\ifx\next\newline%
\gdef\nextthing{\eatone}\fi\fi\fi\egroup\nextthing#1}

{\catcode`\^^M=\active%
\gdef\spacer{\bgroup\catcode`\^^M=\active%
\let^^M=\newline\obeyspaces\lookatspace}}

\def\ref#1{\seeifdefined{#1}\expandafter\csname\one\endcsname\spacer}

\def\cite#1{\expandafter\ifx\csname#1croref\endcsname\relax[??]\else
\csname#1croref\endcsname\fi\spacer}


\def\seeifdefined#1{\expandafter\stripbs\string#1croref*%
\crorefdefining{#1}}

\newif\ifcromessage
\global\cromessagetrue

\def\crorefdefining#1{\ifdefined{\one}{}
{\ifcromessage\global\cromessagefalse%
\message{\spaces\spaces\spaces\spaces\spaces\spaces\spaces}%
\message{<Undefined reference.}%
\message{Please TeX file once more to have accurate cross-references.>}%
\message{\spaces\spaces\spaces\spaces\spaces\spaces\spaces}\fi[??]}}

\def\label#1#2*{\gdef\ctest{#2}%
\xdef\currlabel{\string#1croref}
\expandafter\seeifdefined{#1}%
\ifx\empty\ctest%
\xdef\labelnow{\write\auxfile{\noexpand\def\currlabel{\the\pageno}}}%
\else\xdef\labelnow{\write\auxfile{\noexpand\def\currlabel{#2}}}\fi%
\labelnow}

\def\ifdefined#1#2#3{\expandafter\ifx\csname#1\endcsname\relax%
#3\else#2\fi}




\def\articlecontents{
\vskip20pt\centerline{\bf Table of Contents}\everypar={}\vskip6pt
\bgroup \leftskip=3pc \parindent=-2pc 
\def\item##1{\vskip1sp\indent\hbox to2pc{##1.\hfill}}}

\def\endcontents{\vskip1sp\leftskip=0pt\egroup}

\def\journalcontents{\vfill\eject
\def\currannalsline{\hfill}
\global\titletrue
\vglue3.5pc
\centerline{\tensc\hskip12pt TABLE OF CONTENTS}\everypar={}\vskip30pt
\bgroup \leftskip=34pt \rightskip=-12pt \parindent=-22pt 
  \def\\ {\vskip1sp\noindent}
\def\pagenum##1{\unskip\parfillskip=0pt\dotfill##1\vskip1sp
\parfillskip=0pt plus 1fil\relax}
\def\name##1{{\tensc##1}}}


\institution{}
\onpages{0}{0}
\def\lastpage{???}
\def\thetitle{Title ???}
\def\theauthors{Authors ???}
\def\thereceived{}
\def\therevised{}

\gdef\split{\relaxnext@\ifinany@\let\next\insplit@\else
 \ifmmode\ifinner\def\next{\onlydmatherr@\split}\else
 \let\next\outsplit@\fi\else
 \def\next{\onlydmatherr@\split}\fi\fi\let\eqnu\xspliteqnu\next}

\gdef\align{\relaxnext@\ifingather@\let\next\galign@\else
 \ifmmode\ifinner\def\next{\onlydmatherr@\align}\else
 \let\next\align@\fi\else
 \def\next{\onlydmatherr@\align}\fi\fi\let\eqnu\xspliteqnu\next}

\def\spliteqnu{{\tenrm\sectandeqnum}\relax}

\def\xspliteqnu{\tag\spliteqnu}

\catcode`@=12

\document

\annalsline{April}{1998}
\startingpage{1}     

\comment
\nopagenumbers
\headline{\ifnum\pageno=1\hfil\else \rightheadline\fi}
\def\rightheadline{\hfil\eightit 
The Macdonald conjecture
\quad\eightrm\folio}

\voffset=2\baselineskip
\endcomment


%
%
%
%
%

\def\for{\  \hbox{ for } \ }
\def\if{ \ \hbox{ if } \ }

\def\where{\  \hbox{ where } \ }
\def\and{\  \hbox{ and } \ }
\def\or{\  \hbox{ or } \ }

\def\equal{\buildrel  def \over =}

\def\la{\lambda}
\def\La{\Lambda}

\def\th{\theta}
\def\al{\alpha}
\def\be{\beta}
\def\ga{\gamma}
\def\ep{\epsilon}

\def\de{\delta}

\def\ka{\kappa}
\def\si{\sigma}

\def\Ga{\Gamma}
\def\ze{\zeta}


\def\vph{\varphi}

\def\vep{\varepsilon}

\def\tze{\widetilde{\ze}}

\def\tz{\widetilde z}

\def\hze{\widehat{\zeta}}

\def\C{\bold{C}}
\def\Q{\bold{Q}}

\def\R{\bold{R}}
\def\N{\bold{N}}
\def\Z{\bold{Z}}

\def\one{\bold{1}}

\def\0{\bold{0}}


\def\s{\Cal{S}}

\def\l{\Cal{L}}

\def\p{\Cal{P}}

\def\c{\Cal{C}}

\def\e{\Cal{E}}

\font\germ=eufb10 
\def\goth#1{\hbox{\germ #1}}

\def\EE{\goth{E}}

\def\AA{\goth{A}}
\def\CC{\goth{C}}
\def\GG{\goth{G}}
\def\ZZ{\goth{Z}}
\def\zz{\goth{z}}



\title
{On $q$-analogues of Riemann's zeta}

\shorttitle{On $q$-analogues of Riemann's zeta}

\acknowledgements{
Partially supported by NSF grant DMS--9622829}

\author{ Ivan Cherednik}

\institutions{
Math. Dept, University of North Carolina at Chapel Hill,   
 N.C. 27599-3250
\\ Internet: chered\@math.unc.edu
}


%
%
%
%
%
\vfil

The aim of the paper is to define  $q$-deformations of the
Riemann zeta function and to  extend them 
to the whole complex plane. The construction is directly
related to the recent difference generalization of the 
Harish-Chandra theory of zonal spherical functions [C1,C2,C3].
The Macdonald truncated theta function [M1] replaces $x^s$. 
The analytic continuation is based on the shift operator
technique and a more traditional approach
close to  Riemann's first proof of the functional equation 
(see [E,R]).

We introduce  $q$-deformations using the real and imaginary
integration. The construction depends on the particular choice of
the  path and the integrand much more
than in the classical theory. We come to several different 
functions with different properties. All have meromorphic continuations
to the whole complex plane. The basic ones approach $\ze$ 
(up to a  $\Ga$-factor) as $q\to 1$ for all $s$ thanks 
to the Stirling-Moak formula for $\Ga_q$  [Mo]. 
There are also  equally interesting transitional $q$-zeta functions
converging to $\ze$ for $\Re s>1$ and to  proper combinations of
gamma functions for  $\Re s<1$. They
approximate $\ze$ very well in the critical strip 
when  $q$ is close enough but not too close to $1$. 
Then the tendency
$\Ga_q\to \Ga$   suppresses the Stirling formula and eventually they
go to their $\Ga$-limits.

Our $q$-integrals do not satisfy the celebrated relation 
between $\ze(s)$ and $\ze(1-s)$ 
in  the $q$-setting as well as the Euler product
formula. However certain analytic properties of the $q$-zeta functions
are better than the classical ones. For instance, the imaginary one is
periodic in the imaginary direction and
stable as $\Re s\to \pm \infty$. The number of its zeros can be 
calculated exactly upon the 
symmetrization $s\leftrightarrow 1-s$. Its leading term 
coincides with that from the Riemann estimate of the number
of zeros of $\ze(s)$ in the critical strip. The $q$-zeros
do not belong to the critical line anymore but
their distribution is  far from random.

\vskip 5pt
Let us try to summarize the key results of the paper
in more detail.

1) We construct the  meromorphic continuation of the 
imaginary $q$-zeta functions
(defined in terms of the imaginary $q$-integration) to all $s$
by means of the  shift operator method, mainly Theorem 
\ref\SHIFTWO. The existence of such a continuation is not too
surprising because we can also apply  Cauchy's theorem almost following
Riemann's first proof of the functional equation. The functions with
poles dense everywhere appear in this approach.

2) What is surprising is that the  convergence of the 
meromorphic continuation to the classical zeta (up to a factor)
holds for the $q$-plus-zeta corresponding
to the integration with the kernel $(q^{x^2}+1)^{-1}$, 
but  fails for the  one defined for  $(q^{x^2}-1)^{-1}.$
The latter tends to $\Ga(s)\ze(s)$   for $\Re s>1$ only. 
In the plus case, the limit for negative $\Re s$  is calculated  via 
the shift operator.  We were not able to empoy other methods.

3) Numerical experiments show that the  $q$-plus-zeta
does not have zeros for $\Re s>1/2$ when $0<q<0.96.$ 
Actually it is very surprising 
because the imaginary $q$-zeta functions are periodic in the imaginary
direction and the convergence to the classical zeta (multiplied by
a proper factor) is for relatively  small $\Im s$. Presumably all
zeros of the  $q$-plus-zeta come from  the classical zeros
for $\Re s>0$  and  $a>1.$

4) The imaginary $q$-integral can be exactly calculated for the kernel
$q^{-x^2}.$ It is an important part of the paper
(Section 2) related to  Ramanujan's $q$-integrals and closely
connected with [C1].  
The resulting  function has no zeros for  $\Re s>1/2.$ 
However it explains the above numerical phenomenon for small $q$ only.
Analogous  exact formulas are obtained  for the Jackson summation
instead of the imaginary integration.

5) The most promising result of the paper may be  the definition
of the sharp $q$-zeta functions, which have no classical counterparts.
They are  Jackson sums in terms of 
the  plus-minus kernels discussed above, convergent for
all $s$ except for the poles. Their limits $q\to 1$ are similar
to those for the imaginary integration. The analysis of their
zeros leads to some interesting properties of Riemann's zeta.  

6) Still the main applications are expected upon the
(anti)symmetrization of the integral $q$-zeta functions
with respect to $s\leftrightarrow 1-s$
(Section 6). In the imaginary case, the number of the corresponding
zeros matches well the famous Riemann estimate. The antisymmetrization
of the sharp $q$-zeta may satisfy the Riemann hypothesis.
Conjecturally, all its zeros in a natural horizontal strip are  imaginary.

\vskip 5pt
\vfil

Actually our methods  are of more general nature.
They can be applied to study 
difference Fourier transforms of  functions
of $q^{x^2}$-decay  for 
a quadratic form $x^2$ associated to a root system in $\R^n$
(or in the imaginary variant). Without going into detail, we mention 
that such class of functions is Fourier-invariant in the $q$-theory
in contrast to the Harish-Chandra transform. 
The $q$-Fourier transforms ``at zero'' of
$(q^{-x^2}\pm 1)^{-1},$  simplest 
deformations of the $q^{x^2}, $ are natural
multi-dimensional counterparts of the $q$-integrals studied below.

Note that  $q$-zeta functions  are mainly considered in the paper as    
a good  starting point for the general analytic theory of the
$q$-Fourier transform. We use them to demonstrate the shift operator
method and to show some typical problems. Thanks to
such a concrete  choice  (and the one-dimensional setup)
we are able to illustrate all claims  numerically. There are of
course more specific results  and  applications 
to the classical zeta, which  hopefully will be continued.

The first section contains a certain motivation and the main 
constructions (for the standard minus-kernel).
Section 2 is devoted to the exact formulas for the integrals of the
truncated theta with and without the Gaussian. 
In Section 3, we introduce
the $q$-zeta functions systematically 
using the real and imaginary integrations
and the analytic continuation.
Section 4 contains the approach  via the shift operator.
The next one is a discussion of the $q$-zeros.
The last section is an attempt to establish connections
with the Riemann hypothesis.

The author thanks  P. di Francesco and E. Opdam for useful discussions.
I acknowledge my special indebtedness to David Kazhdan, who greatly
stimulated this work.

%
%
\vfil
\section{Main constructions}

Switching from integers to their $q$-analogues
is a common way of defining $q$-counterparts.
For instance, one can deform  
$\ze(s)=\sum_{n=1}^\infty n^{-s}$ as follows [UN]:
$$\eqalign{
&\ze_q(s)\ =\ \sum_{n=1}^\infty {q^{sn} \over [n]_q^{s}} \for [n]_q\equal
{ 1-q^{n}\over 1-q}.
}
\eqnu
\label\ueno\eqnum*
$$
Here we assume that $0<q<1$ and $\Re s>0$. The numerator 
ensures the convergence. It becomes a constant when
$\ze_q$ is rewritten in terms of $[n]_{q^{-1}}$.
This function  approaches $\ze$ as $q\to 1$ for $\Re s>1.$ We note that
similar $q$-deformations were considered before in connection with 
Carlitz's Bernoulli numbers and $p$-adic $L$-functions (see [Sa]).
Ueno and Nishizawa constructed a meromorphic continuation of $\ze_q$  
to all $s$ with  simple poles
in $\{-\Z_++ 2\pi a i\Z\}$ for $a=-(\log(q))^{-1},\ \Z_+=\{0\le n\in \Z\}$.
The procedure is a certain substitute of
the  functional equation which cannot be saved.
As a rule relations between special
functions which involve the multiplicative property of $x^s$
have no direct $q$-counterparts.  

A  modification of $\ze_q$ tends to $\ze$ for all $s$.
The latter property is important but not sufficient. Once the functional
equation is sacrificed, something else  should be gained.
At least, certain analytic advantages can be expected.
Generally speaking, $q$-functions (and their representations
as sums, products or integrals) are single-valued and
have better analytic properties than their classical counterparts.
It does not seem to happen with $\ze_q(s)$. Actually its definition  is a
transition from the  differential theory to the difference one
because the complex powers are still involved. The latter
are  redefined via $\Ga_q$ in a systematic $q$-theory.

Indeed, substituting $q^{sn}\to q^{sn/2}$ in (\ref\ueno),
the resulting zeta (up to a simple factor) is 
the classical $\ze(s)$ where $x^s$ is replaced by $\sinh^s(x/(2a)).$
By the way, such  numerators are more relevant to associate   $\ze_q(s)$ 
with the quantum $SL_2.$ 
This passage  is well-known in the theory of special functions. 
Instead of the Bessel functions, which are pairwise orthogonal with
respect to the measure $x^s dx$, one gets the hypergeometric and
spherical functions.
The multi-dimensional theory of  $ \sinh(x)^s$ is the cornerstone
of the harmonic analysis on symmetric spaces (see e.g. [He,HO]). 
It is of obvious importance
to incorporate the hypergeometric function into the theory of zeta.
The paper [UN] is a step in this direction.
However the  basic (difference)
hypergeometric function seems more appropriate. 

\vskip 5pt
{\bf A straightforward definition.}
Let us switch from  $ [x]_q^{2k}$ to 
the following
combination of $q$-shifted factorials (see e.g. [A,M1]):
$$\eqalign{
&\de_k(x;q)\ =\ \prod_{j=0}^\infty {(1-q^{j+2x})(1-q^{j-2x})\over
(1-q^{j+k+2x})(1-q^{j+k-2x})}.
}
\eqnu
\label\de\eqnum*
$$

Trying to improve the $ze_q$ due to Ueno - Nishizawa, we put
$$
\eqalign{
&\tze_q^u(s)_c\ =\ \sum_{n=1}^c \de_{-s}(\sqrt{aun};q),\ 
\tze_q^u(s)\ =\ \tze_q^u(s)_\infty\ \for \cr
&q=\exp(-1/a),\ \R\ni a>0,\ u\in\C,\  c=1,2,\ldots .
}
\eqnu
\label\zequ\eqnum*
$$
From now on, $q$ and $a$ will be always connected as in (\ref\zequ).
Since $\de(x)$ is even the sign of $\sqrt{ }$
can be arbitrary.
As we will see,  $\sqrt{a}$ and the ``direction''
$u$ are necessary to ensure the convergence and 
go back to the classical
theory via the following Stirling-Moak formula [Mo,UN].

\proclaim {Lemma} a) Setting  
$ z^k=\exp(k\log z)$ for complex (fixed)
$z, k$ and the usual branch of $\log$,
$$\eqalign{
& \lim_{a\to \infty} (a/4)^{k}\de_k(\sqrt{az};q)\ = 
\ (-z)^k  \cr
& \hbox{if \ either\ } z\not\in \R_+ \hbox{\ or\ } k \in \Z.
}
\eqnu
\label\moak\eqnum*
$$
b) Replacing $az$ by $v$ for any (fixed) $v\in \C$:
$$\eqalign{
& \lim_{a\to \infty} (a)^{2k}\de_k(\sqrt{v};q)\ = 
\ {\Ga(k+2\sqrt{v}) \Ga(k-2\sqrt{v})\over
   \Ga(2\sqrt{v}) \Ga(-2\sqrt{v})}.     
}
\eqnu
\label\galim\eqnum*
$$
\label\MOAK\theoremnum*
\endproclaim

The singularities of the denominator in (\ref\moak) 
go to infinity together with
$a$ for $z\not\in \R_+$ and  disappear in the limit.
When  $k\in \Z$  the claim is obvious since $\de_k$ is a rational
function in this case.
The Moak theorem  also gives a certain estimate of the convergence.
The second formula is well-known.

It is clear why the first formula fails for $z\in \R_+$. The numerator 
of (\ref\de) becomes $0$ for infinitely
many $z$ provided that $k\not\in \Z$. However the limit cannot
be zero. Another reason is that
we switch to (\ref\galim) when $z=v/a$ for fixed
$v$. The $a$-factor in this
formula  is $a^{2k}=o(a^{k})$ for $\Re k<0$. Hence (\ref\moak) 
diverges as $a\to \infty$ for $z$ approaching zero and such $k$.
The latter argument will be applied  to analyze the limiting
behavior of the sharp $q$-zeta functions.

Thus the simplest choice $u=1$ in (\ref\zequ) is impossible, since
it destroys the termwise convergence to the classical sum for $\ze(s)$.
The lemma gives  that
$$ \lim_{a\to \infty} (-4u/a)^s \tze_q^u(s)_c\ = 
\ \sum_{n=1}^c n^{-s}
$$
for any $s\in \C$ (and fixed $c$) if $u\not\in R_+$. However we  must avoid
$u\in -\R_+$ because  
$\lim_{c\to\infty} \tze_q^u(s)_c$ does not exist for such $u.$
Let us examine the imaginary unit $u=i$. 

\proclaim{Theorem}
Let $\tze_q(s)= \tze_q^i(s)$, $\Re s>0$. Then $\tze_q(s)$
exists apart from the set $S(a)$ of the poles of all terms 
in (\ref\zequ). If $\Re s>1$, given a sequence $a_m\to \infty$ and
$k\not\in \cup_m^\infty S(a_m)$, 
$$
\eqalignno{
& \lim_{m\to \infty} (-4i/a_m)^{s}\tze_q(s)\ = 
\  \ze(s)=\sum_{n=1}^\infty n^{-s}.
&\eqnu
\label\limd\eqnum*
}
$$
\label\LIMD\theoremnum*
\endproclaim

We see the first special feature of the $q$-zeta functions.
The convergence is
for $\Re s>0$ in contrast to the usual inequality $\Re s>1$.
However (\ref\limd) holds for  
$\Re s>1$ only.

The function  $\tze_q(s)$ is not real for real $s$. 
The following function is: 
$$
\eqalignno{
&\sum_{n=1}^\infty (\de_{-s}(\sqrt{ain};q)+
\de_{-s}((\sqrt{ain})^*;q))/2.
&\eqnu
\label\tzere\eqnum*
}
$$
Here and further $z^*$ is the complex conjugation.
The limit of this function as $a\to \infty$ remains the same.

The analytic behavior of  $\tze_q(s)$ is totally singular. The 
set $S(a)$ is dense everywhere in $\C$ unless  $2\pi a \in \Q$.
If $2\pi a=n/d$ for coprime integers $n,d$, then $S(a)$ belongs
to the union $\c_d$ of the translations of the
diagonal cross $\c=\{x^2\in i\R\}$ by $\{(l/d)i,\ 
 l\in \Z\}$.
Choosing $a_m=m^2 a$, we have the embeddings
 $S(a_m)\subset S(a)$ and  can apply the
theorem at least when  $s\not\in \c_d$. In general, 
the sum for  $\tze_q(s)$ is  pointwise convergent in spite of
infinitely many poles in any neighborhood 
since the residues are getting
smaller faster than the distances to any fixed $s\not\in S(a)$.

It seems that ``totally singular'' functions  are
inevitable in the difference calculus. 
However we prefer to
avoid them in the basic definitions. The next step is to 
remove as many poles as possible turning from sums to
integrals. We note that  $\tze_q(s)$ will be
directly related to  $\ZZ^\la_{+q}(k|d)$ from
(\ref\zzbul) playing important role in the analytic
continuation of the imaginary $q$-zeta functions.

\vskip 5pt  {\bf Real integration.}
Actually the $\de$-function and its multi-dimensional
generalizations due to Macdonald are important because they
make selfadjoint remarkable difference operators which have
many applications (see [M2,C3]). Hence it is 
more logical to use them as measures, which 
is exactly what we are going to do.

From now on we will switch from  $s$ to the ``root multiplicity''
$k=s-1/2$, standard in the harmonic analysis on symmetric spaces 
 and in the corresponding theories of orthogonal polynomials. This choice
is convenient in the theory of  $\ze$  as well: 
the critical line $\Re s=1/2$ coincides with the $k$-imaginary axis.

First, we will express the $\Ga$-function 
 and $\ze\Ga$ in terms of the Gaussian
$\exp(-x^2)$ modifying a little the classical representations:
$$
\eqalignno{
& \Ga(k+1/2)\ =\ 2\int_{0}^\infty x^{2k}\exp(-x^2) dx, \ \Re k>1/2, \cr
& \ZZ(k)\equal\ze(k+1/2) \Ga(k+1/2)\ =
\ 2\int_{0}^\infty x^{2k}(\exp(x^2)-1)^{-1} dx.
&\eqnu
\label\zega\eqnum*
}
$$
Since the first formula has a nice  $q$-generalization with 
$\de_k(x;q)$ instead of
$x^{2k}$ (see [C1] and  the next section), we can try to use the second
to deform $\ZZ(k)$. Let
$$
\eqalignno{
&\ZZ^{\pm}_q(k) \equal \ \int_{0}^{\infty\pm\ep i} 
(\exp(x^2/a)-1)^{-1} \de_k(x;q) dx.
&\eqnu
\label\ZZq\eqnum*
}
$$
The contour  of integration $C^\ep_{\pm}$
for $\ZZ^{\pm}$ goes from $0$ to $\pm \ep i$
and then to $\infty \pm \ep i$ for  fixed $\ep>0$, which is assumed
to be sufficiently small. 
Note that we cannot integrate
over $\R_+$ for real $k$ because $\de_k$ has  infinitely many
poles there. In contrast to (\ref\zega), there is no problem 
with  the point $x=0$. Indeed, $\de_k$ has a zero of the second 
order at $0$ which annihilates the pole of the same order coming from
 $(\exp(x^2/a)-1)^{-1}$. To be more exact, we need to exclude $k=0$
in this argument. The functions   $\ZZ_q^{\pm}$ 
are regular  apart from 
$$
\eqalignno{
& K^{\ep}_\pm(a)\ = \ 
\{2C^\ep_\pm -\Z_+ +(2\pi a i)\Z\}\cup \{-2C^\ep_\pm -\Z_+ +(2\pi a i)\Z\}.
&\eqnu
\label\srepm\eqnum*
}
$$
They are obviously $2\pi a i$-periodic. However it is of no use, since
the connected components of  $\C\setminus K^{\ep}_\pm(a)$
have the same periodicity. In this definition, we loose the reality
of $\zeta$ on the real axis. The following functions have this important 
property:
$$
\eqalignno{
& \ZZ_q^{re}(k)\equal {\ZZ_q^-(k)+\ZZ_q^+(k)\over 2} =
{1\over2}\int_{-\infty-\ep i}^{+\infty-\ep i}{ \de_k(x;q)\  dx\over
\exp(x^2/a)-1},
&\eqnu\cr
\label\zzre\eqnum*
&\ZZ^{\sharp}_q(k) \equal \ {\ZZ_q^-(k)- \ZZ_q^+(k)\over2i} =
{1\over 2i}\int_{\infty+\ep i}^{\infty-\ep i } {\de_k(x;q)\ dx
\over \exp(x^2/a)-1}.
&\eqnu
\label\zzsharp\eqnum*
 }
$$
The second integration is with respect to the path $\{C^\ep_- - C^\ep_+\}$.
The corresponding function  is well-defined  in the horizontal
strip $\{-2\ep<\Im k<+2\ep\}$ and has no singularities 
for $\Re k>0$.
The natural domain of the first function is 
$K^{re}\equal\C\setminus\{K^{\ep}_-\cup K^{\ep}_+\}$.

\proclaim{Theorem}
a) The functions $\ZZ_q^{\pm}$ coincide when
$$
\eqalignno{
&\Im k +2\pi a\Z\not\in [-2\ep,+2\ep].
&\eqnu
\label\imkep\eqnum*
}
$$
The functions  $\ZZ_q^{re},\ZZ_q^{\sharp}$ have  meromorphic continuations 
from $\R_+$to $\C$. The first  has no poles in the strip
$$
\eqalignno{
&K^{\sharp}(a)\equal\ \{k,\ -2\ep_a<\Im k<2\ep_a\} \for
\ep_a\equal\min\{\sqrt{\pi a},\pi a/2\}.
&\eqnu
\label\epa\eqnum*
}
$$
The set of poles of the continuation of
 $\ZZ_q^{\sharp}$ in this strip  is $\{-\Z_+/2\}$.
It also has  simple zeros at $k=1,2,\ldots\ $.

b) Let   $\ZZ(k)=\ze(k+1/2)\Ga(k+1/2)$ (as above)
and $\Re k>1/2$. Then
$$
\eqalignno{
&\lim_{a\to \infty} (a/4)^{k-1/2}\exp(\pm\pi ik)\ZZ_q^{\pm}(k)\ =\ 
\ZZ(k).
&\eqnu\cr
\label\limzpm\eqnum*
&\lim_{a\to \infty} (a/4)^{k-1/2}\ZZ_q^{\sharp}(k)\ =\ 
\sin(\pi k)\ZZ(k).
&\eqnu
\label\limsharp\eqnum*
}
$$
For any $k\in\C$,
$$
\eqalignno{
&\lim_{a\to \infty} (a/4)^{k-1/2}\ZZ_q^{re}(k)\ =\ 
\cos(\pi k)\ZZ(k).
&\eqnu
\label\limzre\eqnum*
}
$$
\label\LIMRE\theoremnum*
\endproclaim

The calculation of the limits of $\ZZ_q^{\pm}$ is based
directly on (\ref\zega) and Lemma \ref\MOAK; (\ref\limsharp)
is a formal corollary of the (\ref\limzpm).
 
The limit of  $\ZZ_q^{re}(k)$ is calculated using the following
well-known integral representation for the zeta (see e.g. [E]): 
$$
\eqalignno{
&(1/2)\oint_{\infty-\ep i}^{\infty+\ep i}
{(-z)^k dz\over i\sqrt{-z}(e^z-1)}\ =\ -\sin(\pi s)\ZZ(k)=
\cos(\pi k)\ZZ(k).
&\eqnu
\label\intzz\eqnum*
}
$$
The path of integration 
begins at $z=-\ep i+\infty$, moves
to the left down the postive real axis till $-\ep i$, then circles the origin
and returns up the positive real axis to $\ep i+\infty$
(for small $\ep$). Note the difference from 
the integral $\int_{\infty+\ep i}^{\infty-\ep i}$ above, 
where the orientation is different and the path goes through the origin.
The left-hand side multiplied by $(a/4)^{k-1/2}$
is exactly the integral (\ref\zzre) after
the substitution: $z=x^2/a$. The factor $i\sqrt{-z}$
is $(a/2)dz/dx$, $\sqrt{x}$ is the standard branch positive on $R_+$. 

In the paper, we prefer to 
replace $(\exp(x^2/a)-1)^{-1}$  by $(\exp(x^2/a)+1)^{-1},$ 
which improves the analytic properties of the $q$-zeta functions.
For instance, the left-hand side of (\ref\limsharp) diverges
for $\Re k<1/2,$ whereas its ``plus-counterpart''
approaches the zeta (up to a proper factor)
for all $k.$
However the asymptotic behavior of $\ZZ_q^{\sharp}(k)$ 
for $\Re k<1/2$ is very interesting:
$$
\eqalignno{
&\lim_{a\to \infty} a^{2k-1}\ZZ_q^{\sharp}(k)\ =\ 
\tan(\pi k)\Ga(k)^2 \for  \Re k<1/2, 
&\eqnu\cr
\label\limshn\eqnum*
&(a/4)^{k-1/2}\ZZ_q^{\sharp}(k)= 
(4a)^{1/2-k}\tan(\pi k)\Ga(k)^2 +\sin(\pi k)\ZZ(k)+o(1) 
&\eqnu
\label\limsho\eqnum*
}
$$
for $\Re k=1/2$, where
the error $o(1)$ tends to $0$ as $a\to \infty$. In fact the latter
formula holds for all $k$. Indeed, the term $ \sin(\pi k)\ZZ(k)$  
is suppressed in (\ref\limsho) for $\Re k>1/2$ 
and vice versa for $\Re k<1/2$ 
because of the difference in the   $a$-factors. Moreover,
the error is of order $\tan(\pi k)\Ga(k)^2$ in the strip
$-1/2+\ep<\Re k<1/2-\ep$ for any small $\ep>0$, 
when $a$  are big  but not too big. Since $\Ga(k)^2$ quickly
approaches $0$ as $\Im k$ increases, the function 
$(a/4)^{k-1/2}\ZZ_q^{\sharp}(k)$ first approximates $\sin(\pi k)\ZZ(k)$
with high accuracy and then slowly switches to the 
the leading term $(4a)^{1/2-k}\tan(\pi k)\Ga(k)^2$ in the critical strip.
Numerically, the latter (dominant) regime is very difficult 
to reach for big $\Im k$, even for  $k \sim 10i$. 

\comment
To check these formulas and evaluate  the error we use the following
variant of (\ref\intzz):
$$
\eqalignno{
&\oint_{\infty-\ep i}^{\infty+\ep i}
{(-z)^{2k} dz \over (e^{z^2}-1)}\ =\ i\sin(2\pi k)\ZZ(k) \ 
(\hbox{\ all\ } k),
&\eqnu
\label\intzz\eqnum*
}
$$
substituting $x=z\sqrt{a}$ in (\ref\zzsharp) and calculating
the difference $\oint-\int$ via the residues.
\endcomment

The main advantage of  $\ZZ_q^{\sharp}(k)$ is the following 
representation  (Cauchy's theorem),
which has no counterpart in the classical theory:
$$
\eqalignno{
&\ZZ_q^{\sharp}(k)\ =\ 
-{a\pi\over 2}\prod_{j=0}^\infty
{(1-e^{-(j+k)/a})(1-e^{-(j-k)/a})\over
 (1-e^{-(j+2k)/a})(1-e^{-(j+1)/a})}\times\cr
&\sum_{j=0}^\infty\prod_{l=1}^j
{(1-e^{-(l+2k-1)/a})(1-e^{(l+k)/a})\over
 (1-e^{-(l+k-1)/a})(1-e^{l/a})}{1\over \exp({(k+j)^2\over 4a})-1}.
&\eqnu\cr
\label\caumin\eqnum*
}
$$
Thanks to the factor $(\exp((k+l)^2/(4a))-1)^{-1}\sim \exp(-l^2/(4a))$,
it converges for all $k$ ($0<q<1$) except for  the  singularities,
and very quickly.  For instance, we
can caluclate $\ZZ_q^{\sharp}(k)$ 
numerically for $a\sim 10000$ (i.e. when $q\sim 0.9999$), and $\Im k\sim 100$.
For integral $k$, the range is almost unlimited.
The numerical simulation of  the $\ZZ_q^{re}(k)$ is approximately
$10 \cdots 100$ times worse.

In a sense, the series (\ref\caumin) is complimentary to
$\tze_q(s)$ (Theorem \ref\LIMD). They appear together in the following
construction.

\vskip 5pt {\bf Imaginary integration.}
The integration in the real and imaginary directions are practically
interchangable in the classical theory of zeta function
(the substitution $x\to ix$ multiplies $x^s$ by a constant).
In the harmonic analysis, the difference is fundamental.
The imaginary case (integrating with $\sin(x)^{2k}$)
is much simpler than the real one (for $\sinh(x)^{2k}$ as the kernel).
The real and imaginary integrations are closer to each other
for $\de_k(x;q)$. However the latter has certain  analytic
advantages. 

For $\Re k>0$, we set
$$
\eqalignno{
&\ZZ^{im}_q(k) \equal \ (-i)\int_{0}^{\infty i} 
(\exp(-x^2/a)-1)^{-1} \de_k(x;q) dx.
&\eqnu
\label\zzim\eqnum*
}
$$
The complex conjugation is  denoted by  $*$.

\proclaim{Theorem}
a) The function $\ZZ^{im}_q(k)$ is analytic for $\Re k>0$, and has
a $2\pi a i$-periodic 
meromorphic continuation $\ZZ^{an}_q(k)$ to all $k\in \C$ 
with the following set of
poles (simple for generic $a$):
$$
\eqalignno{
&\p \ =\ \{-\Z_+/2 +\pi a i\Z\}\cup \La, \for
&\eqnu\cr
\label\zzpol\eqnum*
\La = \{-2&\sqrt{2\pi a i \N}-\Z_+ +2\pi a i\Z\}\cup
\{-2\sqrt{2\pi a i \N}^{\,*}-\Z_+ +2\pi a i\Z\}.
}
$$
The poles of the ($2\pi a i$-periodic) 
meromorphic continuation $\widetilde{\ZZ}^{an}_q$
of the function
$$
\eqalignno{
&\widetilde{\ZZ}^{im}_q(k)=\ZZ^{im}_q(k) 
\prod_{j=0}^\infty {(1-q^{j+2k})(1-q^{j+1})\over(1-q^{j+k})(1-q^{j+k+1})}.
&\eqnu
\label\zzstab\eqnum* 
}
$$
belong to the set  $\widetilde{\p}=\La\cup\{-\Z_+ +2\pi a i\Z\}$
(they are double for $\Re(z)\in -\N= -1-\Z_+$).

b)If $\Re k> 1/2$, then 
$$
\eqalignno{
&\lim_{a\to\infty}(a/4)^{k-1/2}\ZZ^{im}_q(k)\ =\ 
\ZZ(k)=\Ga(k+1/2)\ze(k+1/2),
&\eqnu\cr
\label\zziml\eqnum* 
&\lim_{a\to\infty} a^{k-1/2}\widetilde{\ZZ}_q^{im}(k) = 
\sqrt{\pi}\,\Ga(k+1)\ze(k+1/2).
&\eqnu
\label\zzstabl\eqnum*.
}
$$
For all $k$,
$$
\eqalignno{
({a\over 4})^{k-1/2}\ZZ^{an}_q(k)\ \cong\
&(4a)^{1/2-k}\tan(\pi k)\Ga(k)^2 +\ZZ(k),
&\eqnu
\label\limnegim\eqnum*
}
$$
as  $a\to\infty$.

c) The function $\ZZ^{im}_q(k)$ tends to a nonzero limit  as 
 $\Re k\to \infty$.
Given $\ka\in-\R_+$ such that $\{\ka+i\R+\Z\}\cap\widetilde{\p}=\emptyset$,
the limit
$$
\eqalignno{
& \psi(\nu;\ka)\ =\ 
\lim_{\Z\ni m\to \infty}\widetilde{\ZZ}_q^{an}(\ka-m+i\nu)
&\eqnu
\label\psinu\eqnum* 
}
$$
exists and is a continuous $2\pi a$-periodic function of $\nu\in \R$.  
\label\ZZIM\theoremnum*
\endproclaim

The term $\ZZ$ in (\ref\limnegim) can be skipped for $\Re k$ 
strictly less than $1/2$ due to the $a$-factor. However,  the 
formula with this term gives more than was stated in the theorem.
In the critical strip $\{-1/2<\Re k<1/2\}$, the right-hand side 
approximates  well  the left-hand side  even
when $a$  is not too big . For instance, 
 $(a/4)^{k-1/2}\ZZ_q^{an}$ is getting very close  to $\ZZ$ 
when $\Im k\sim 10$ or more  to make $\Ga(k)^2$ small enough. 
It is the same phenomenon (with the same $\Ga$-limit) 
as in (\ref\limsho).

The set of admissible $\ka$ in c) is dense everywhere in $-\R_+$
because  so is $\widetilde{\p}+i\R+\Z\subset\C$
(for generic $a$). However the function $\psi(\nu;\ka)$ is
well-defined. It is  related to the existence of
the totally singular  zeta $\tze_q$ considered above. In fact,  $\psi$
is a limit of a certain sum similar to (\ref\zequ).

Thanks to c),  the number of zeros of  $\ZZ^{an}_q(k)$ in
the vertical strip
$$\{\ka-m-1<\Re k<\ka-m\} \hbox{\ modulo\ }  2\pi a i\Z$$
 equals  the number of poles of $\widetilde{\ZZ}^{im}_q(k)$
(counted with multiplicities)
in this strip for sufficiently big $m$.
It is likely that  even
the positions of individual  zeros (modulo $\Z+2\pi a i\Z$)
approach certain limits  when $m\to\infty$. This can be seen
numerically. 

The number of zeros of  $\ZZ^{an}_q$
in the  half-planes $\Re k>C$ for any $C\in \R$ is finite, 
which readily results from a). It is not difficult to calculate
$\psi(\nu;\ka)$ numerically for $a\sim 10$ or less.
The change of  $\arg(\psi)/(2\pi)$
from $\nu=0$ through $\nu=2\pi a$  combined with the periodicity
of  $\ZZ^{an}_q$ can be used  to calculate the exact
number of its zeros  modulo the period $2\pi a i$
in the half-planes  $\Re k>C.$

We can ensure the  convergence to Riemann's zeta
for all $k$ switching to the  ``plus-zeta functions'' defined
for $(\exp(-x^2/a)+1)^{-1}.$ 
In the classical theory, this substitution 
adds a simple factor
to $\ze$, terminating the pole at $s=k+1/2=1$ and  
improving the  analytic properties.
The zeros for $\Re k\neq 1/2$ remain the same. Setting
$$
\eqalignno{
&\ZZ^{im}_{+q}(k) \equal \ (-i)\int_{0}^{\infty i} 
(\exp(-x^2/a)+1)^{-1} \de_k(x;q) dx \for \Re k>0,
&\eqnu\cr
\label\zzimpl\eqnum*
&(a/4)^{k-1/2}\ZZ^{im}_{+q}(k)\ \to \  
\ZZ_+(k)\equal (1-2^{1/2-k})\ZZ(k).
}
$$
The limit holds for all $k\in \C$ except for the poles
if $\ZZ^{im}_{+q}$ is replaced by the meromorphic continuation
$\ZZ^{an}_{+q}(k).$

\vskip 5pt  {\bf Locating the zeros.}
The numerical experiments demonstrate  that  $\ZZ^{im}_q$
has only two zeros ($\mod 2\pi a i$) in the right half-plane $\Re k>\ep$
at least for $\ep=0.05,\ a\le 25 (q<0.96)$. These  zeros are connected
to each other: 
$$ z_0=\xi+\eta i,\ \  z_1=\xi+(2\pi a-\eta)i,\ \eta\le \pi a.$$
In a sense, these zeros appear because of
the poles of  $\ZZ_q^{im}$ at $k=0,\pi a i$.
They vanish for $\ZZ^{im}_{+q}(k)$. 

The $z_0$
slowly approaches $k=1/2$, the pole of $\ze(k+1/2)$:
$$
\eqalignno{
&z_0(a=3) = 0.0215+ 1.2746i, \ z_0(100)\ =
0.4419+ 0.7216i\cr
&z_0(1000)\, =0.4741+0.5560i, \ z_0(2000)=
0.4788+0.5208i.
&\eqnu
\label\strzim\eqnum*
 }
$$

Switching to  $\ZZ^{an}_{+q},$ we hope that its zeros  in  the strips
$\Re k>-\ep,\ \ep>0$
are always deformations of the zeros of  $\ze(k+1/2)$ for sufficiently 
big $a>a(\ep)$. This might be true for $\ZZ^{an}_{q}$ as well with a
reservation about  the ``low''
zeros going towards  $k=1/2$
(similar to  (\ref\strzim)) and their reflections. 

Let us mention without going into detail, that
the following ``theta''-deformation
$$
\eqalignno{
&\ZZ^{\th}_{q}(k) \equal \ (-i)\int_{0}^{\infty i} 
(\sum_{n=1}^{\infty}\exp(-(nx)^2/a)) \de_k(x;q) dx \for \Re k>0
&\eqnu\cr
\label\thimzl\eqnum*
}
$$
also has only low zeros in the same range of $a$.
Here the functional equation for the theta is used to
integrate near $0$. The analytic continuation of $\ZZ^{\th}_q$
to the left half-plane is not known.

We note that the roles of the left and right half-planes are not symmetric
in the $q$-theory. For instance, the 
zeros of  $\ZZ^{an}_{+q}$  mainly appear in the left one.
There is an explanation for this.

Given a zero $k=z\in i\R$ of $\ze(k+1/2)$, the linear approximations 
$\tz_+(a)$ 
for the zero $z_+(a)$ of $\ZZ^{an}_{+q}$ associated to a given
zero $k=z$ of $\ze(k+1/2)$ can be calculated as follows:
$$
\eqalignno{
& \tz_+(a)=z(1+{ 4(z+1/2)\ze_+(z+3/2)-(z-1)\ze_+(z-1/2)
\over 12 a\ze'(z+1/2)(1-2^{1/2-z}) }) \cr
&\for \ze'(s)=\partial \ze(s)/\partial s,\ \ze_+(s)=(1-2^{1-s})\ze(s).
&\eqnu
\label\tzappl\eqnum* 
}
$$

The formula results  from (\ref\zzimpl).
However we did not check the necessary estimates, so it is 
not completely justified at the moment. 
The first $\tz_+(a)$ which we found in the right half-plane 
corresponds to $z= 1977.2714i$, well beyond the range where we
can calculate  $\ZZ^{an}_{+q}$ numerically.
After this point, the number of linear deformations
in the right half-plane grows but not quickly. Generally speaking,
they appear when the distances between the corresponding
consecutive zeros of $\ze$
are getting too close, i.e. exactly when the linear approximation could not be
reliable (numerically).

The same strong tendency 
holds for the ``real'' zeta $\ZZ^{re}_q$. However now the
  deformations of the zeros prefer
 the left half-plane. 
The formula for the linear approximation  $ \tz(a)$ of the zero
$z(a)$ of  $\ZZ^{re}_q(k)$ corresponding to a given root $z$ of $\ze(k+1/2)$
reads: 
$$
\eqalignno{
& \tz(a)=z(1-{ 4(z+1/2)\ze(z+3/2)-(z-1)\ze(z-1/2)
\over 12 a\ze'(z+1/2)}).
&\eqnu
\label\tzan\eqnum* 
}
$$
Note the opposite sign of the  term after $1$
and the missing factor $(1-2^{1-s})$. The formula  for the plus-counterpart
 $\ZZ^{re}_{+q}(k)$ of  $\ZZ^{re}_q(k)$ 
is   (\ref\tzappl) with the minus after $1$. The 
tendency to move to the right is  practically the same
for $\ZZ^{re}_{+q}(k).$ The 
factor $(1-2^{1-s})$  does not contribute too much.

Formula (\ref\tzan) gives a very
good approximation for the zeros $z^{\sharp}(a)$ of
$\ZZ^{\sharp}_q(k)$ as well. Of course $a$ must be big but not enormously
big, because eventually  $\ZZ^{\sharp}_q(k)$ will go to its $\Ga$-limit, 
which has no zeros in the critical strip.
For instance, $\tz(100)= 0.7091 + 14.1326i$ for 
$z= 14.1347i$, the first zero of $\ze(k+1/2)$. The exact zero
$z^\sharp(100)$ equals $0.7041 + 14.0779i$. Enlarging $a$:
$$ \tz(1000)=0.0709 + 14.1345i,\  \   z^\sharp(1000)= 0.0709 + 14.1340i. $$
 

The first zero $z$ of $\ze(k+1/2)$ with the ``real''
linear approximation $\tz$
in the left half-plane (which we found) is $z=1267.5706i$. The next one
is $1379.6833i$. These zeros are  beyond the range we are able to analyze
numerically, although the calculations with $\ZZ^{\sharp}_q$ 
are $10\cdots 100$ times faster than those  for $\ZZ^{re}_q,\ \ZZ^{im}_q$.

The $a$-evolution of the ${\sharp}$-zeros is very  nontrivial.
It is not too difficult to trace the evolution of the first ones
back in $a$ (at least, through $a\sim 1$). For instance, if $z= 14.1347i$
$$
\eqalignno{
&z^{\sharp}(35)=1.92 + 13.71i,\
z^{\sharp}(11)\ = 4.25+11.36i,&\eqnu\cr
\label\strsmall\eqnum*
&z^{\sharp}(5)\ \, =0+8.72i,\ \ \ \ \ \
z^{\sharp}(2)\ \  =-2.98+5.01i.
}
$$

Concerning big $a$, we were not able to reach the moment when
the zero under consideration leaves the vicinity of
$z= 14.1347i$. The linear approximation
works reasonably
well at least untill $a=50000$. One  may expect $z^{\sharp}$ to go 
to the vertical line $\Re k=1/2.$  Employing the formula
(\ref\limsho), $z^{\sharp} \simeq 0.0015+14.1320i$
when $\log(a)= 122.25,$ so it is still close to the classical zero
for such $a$. When  $\log(a)= 1020.34,$ approximately
$z^{\sharp}\simeq 0.4361+ 14.1334i$. This 
calculation is of course qualitative (as well as(\ref\limsho)), 
but presumably  demonstrates  the tendency.

There are several (more than $1$) examples
of ``low'' zeros of $\ZZ^{\sharp}_q$ 
traveling from the left half-plane  towards 
the natural ``attractor'' $k=1/2.$
The evolution of the first of them (a counterpart of
the zero from (\ref\strzim)) is as follows:
$$
\eqalignno{
a=5: &z^{\sharp}= 0.037 + 0.735i,\ a=\ \ 2733: z^{\sharp}=
 0.409 + 0.477i,&\eqnu\cr
\label\strzre\eqnum*
 a=\ \ 4500: &z^{\sharp}=0.418 + 0.460i,\
a=\ 15187: z^{\sharp}=0.436 + 0.424i,\cr 
a=\ 76882: &z^{\sharp}=0.452 + 0.383i,\    
a=172984: z^{\sharp}=0.459 + 0.365i.
}
$$  
The convergence is very slow. It explains to a certain
extent why the deformations of the zeros of $\ze$
stay put. The speed of decay at big $a$ is expected
to be much  smaller for the ``high''
zeros.

Summarizing, the distribution of the zeros of
the imaginary $q$-zeta functions looks surprisingly regular.
As to the  $\ZZ^{\sharp}_q$ and its
plus-counterpart,
the main objects of our computer simulation, the $q$-zeros are 
also far from random
in the horizontal strip $ \{-2\ep_a<\Im k<2\ep_a\}$ for
$\ep_a=\sqrt{\pi a}$ (see (\ref\epa)). 

There must be a reason for this. Numerically,   
a little decrease of the accuracy  
in the formula  
for $\ZZ^{\sharp}_q$ (say, $\sim (5\cdots 10)\sqrt{a}$ terms in the range
$k= 14i\cdots 50i$ instead of the 
necessary $\sim (15\cdots 20)\sqrt{a} $ for $a\sim 500\cdots 1000$) 
and one can see how irregular  the  zeros could  be.

%
%
%
\vfil
\section{Exact formulas}

First we  will drop the factor $\exp((k+l)^2/(4a))-1)$ 
in the  formula (\ref\caumin) for $\ZZ_q^{\sharp}(k)$ (see
(\ref\zzsharp)) and calculate
the resulting sum. It is a  variant of
the celebrated constant term conjecture  [AI,M1,M2,C3].
The latter  corresponds to the imaginary 
integration over the period.
The following  theorem  readily results from the $q$-binomial theorem.
Let $\ep>0$,
$$
\eqalignno{
&\CC_q^{\sharp}(k)\equal
{1\over 2i}\int_{\infty+\ep i}^{\infty-\ep i} \de_k(x;q)\  dx\ =
&\eqnu\cr
\label\ccint\eqnum*
&-{a\pi\over 2}\prod_{j=0}^\infty
{(1-e^{-(j+k)/a})(1-e^{-(j-k)/a})\over
 (1-e^{-(j+2k)/a})(1-e^{-(j+1)/a})}\times\cr
&\sum_{j=0}^\infty \prod_{l=1}^j
{(1-e^{-(l+2k-1)/a})(1-e^{(l+k)/a})\over
 (1-e^{-(l+k-1)/a})(1-e^{l/a})}.
&\eqnu
\label\ccdef\eqnum*
}
$$
The last series is convergent for 
all $k$ except for $k\in -1/2-\Z_+ +2\pi a i\Z$.
The coincidence of the integral and the series holds provided
that $k/2$ is inside the contour  of integration.
We will always use the series if $k$ is outside.
 
\proclaim{Theorem}
For all $k$ such that $\Re k <0,\ \ -{2\over a\pi}\CC_q^{\sharp}(k)=$
$$
\eqalignno{
2\prod_{j=0}^\infty
{(1-e^{-(j+k)/a})(1-e^{-(j-k)/a})\over
 (1-e^{-(j+2k)/a})(1-e^{-(j+1)/a})}
\sum_{j=0}^\infty e^{kj/a}\prod_{l=1}^j
&{1-e^{-(l+2k)/a}\over
 1-e^{-l/a}}\ =\ \cr
2\prod_{j=0}^\infty
{(1-e^{-(j+k)/a})(1-e^{-(j+k+1)/a})\over
 (1-e^{-(j+2k)/a})(1-e^{-(j+1)/a})}&.
&\eqnu
\label\ccprod\eqnum*
}
$$
The limit of the latter expression  as $a\to \infty$ is
${2\, \Ga(2k)\over \Ga(k)\Ga(k+1)}.$
\label\CT\theoremnum*
\endproclaim

\vskip 5pt  {\bf Integrating the Gaussians.}
Now let us integrate the Gaussians $q^{\pm x^2}$ with respect to $\de_k$.
See [M1,O,J] about the classical (i.e. differential)
multi-dimensional generalizations.
The formula for  the imaginary integral  can be 
deduced from [C1] (in the case of $A_1$). Let us discuss it first.

\proclaim{Theorem} Provided that $\Re k>0$,
$$
\eqalignno{
&\GG_q^{im}(k)\equal
{1\over 2i}\int_{-\infty i}^{+\infty i}e^{x^2/a}\de_k(x;q)\  dx =
\sqrt{a\pi}\prod_{j=0}^\infty
{1-e^{-(j+k)/a}\over
 1-e^{-(j+2k)/a}},
&\eqnu\cr
\label\ggim\eqnum*
&{a}^{k-1/2}\GG_q^{im}(k)\ \to \ \sqrt{\pi}\, {\Ga(2k)\over \Ga(k)}
\hbox{ \ as \ } a\to\infty. &\eqnu
\label\ggimlim\eqnum*
}
$$
\label\GGIM\theoremnum*
\endproclaim

\vfil
Actually the theorem can be considered as 
a variant of the Ramanujan integrals, namely  (3.12) and (3.13)
from Askey's paper in [A]. They  correspond  to the
two  special cases $k=\pm \infty$, when $\de_k$ has either no
denominator or numerator, but involve more parameters.
Richard Askey writes, ``try to imagine what someone with Ramanujan's
ability could do now with the computer algebra systems ...''.
Of course the computers are useful for algebraic $q$-calculations,
but still supporting the $q$-analytic theory (meromorphic 
continuations, zeros etc.) seems a more promising application.
The theory is not mature yet and we need numerical
experiments. However, first, computers  must become
at least $100$ times faster.

  
We note that the integral from (\ref\ggim) does not give
the analytic continuation of $\GG_q^{im}(k)$ for   $\Re k<0$. 
There is a simple way to 
extend the theorem  to $\Re k>-1/2$ replacing $\de_k$
by its asymmetric variant
$$
\eqalign{
&\mu_k(x;q)\ =\ \prod_{j=0}^\infty {(1-q^{j+2x})(1-q^{j+1-2x})\over
(1-q^{j+k+2x})(1-q^{j+k+1-2x})}.
}
\eqnu
\label\mu\eqnum*
$$
The $\mu$-measure has many applications and appears
naturally in the proof of the Macdonald constant term conjecture
(formulated in terms of $\de$ only). Actually it is not surprising
because the symmetry $x\to 1/2-x$ of $\mu$ is quite common
for many $q$-functions (including the theta functions).
The connection with the $\de$-measure  is through the symmetrization:
$$
\eqalignno{
&\int_{-\infty i}^{+\infty i}f(x^2)\mu_k\  dx\  =\
{1+e^{-k/a}\over 2}\int_{-\infty i}^{+\infty i}f(x^2)\de_k\  dx
\hbox{\ for\ any\ } f.
&\eqnu\cr
\label\mugg\eqnum*
}
$$
The function 
$$
\eqalignno{
&\GG_q^{\mu}(k)\equal
-{i\over 1+e^{-k/a}}\int_{1/4-\infty i}^{1/4+\infty i}e^{x^2/a}
\mu_k(x;q)\  dx
&\eqnu\cr
\label\immu\eqnum*
}
$$
coincides with $\GG_q^{im}(k)$ for $\Re k >0$ and gives its 
meromorphic continuation to $\Re k>-1/2$ with the poles
$\{\pi a i+2\pi a i\Z\}$ coming from the factor before the integral.
Respectively, (\ref\ggim) holds for such $k$ if  $\GG_q^{im}(k)$
is replaced by  $\GG_q^{\mu}(k)$. The passage $\de\to \mu$ will be
also used later for the  zeta.

{\bf Real integration.} A more interesting exact formula is for
$$
\eqalignno{
&\GG_q^{\sharp}(k)\equal
{1\over 2i}\int_{\infty+\ep i}^{\infty-\ep i} e^{-x^2/a} \de_k(x;q)\  dx\ =
&\eqnu\cr
\label\ggint\eqnum*
&-{a\pi\over 2}\prod_{j=0}^\infty
{(1-e^{-(j+k)/a})(1-e^{-(j-k)/a})\over
 (1-e^{-(j+2k)/a})(1-e^{-(j+1)/a})}\times\cr
&\sum_{j=0}^\infty e^{-{(k+j)^2\over 4a}} \prod_{l=1}^j
{(1-e^{-(l+2k-1)/a})(1-e^{(l+k)/a})\over
 (1-e^{-(l+k-1)/a})(1-e^{l/a})}.
&\eqnu
\label\ggdef\eqnum*
}
$$
Here the convention is the same as for $\CC_q^{\sharp}(k)$:
we use the integral if and only if it coincides with series, i.e.
when $k/2$ is inside the contour of integration. For such $k$,
(\ref\ggdef) is nothing else but an application of the
Cauchy theorem.

\proclaim {Theorem} Let $Q=\exp(-8a\pi^2)$. Then for all $k$
except for the poles
$$
\eqalignno{
\GG_q^{\sharp}(k)\ &=\ \phi(k)\sqrt{\pi a}\, \exp({k(k+1)\over 2a})
\prod_{j=0}^\infty
{1-e^{-(j+k)/a}\over
 1-e^{-(j+2k)/a}} \for\cr
\phi(k)\ &=\ \sin(\pi k)\prod_{j=1}^\infty {1+Q^{j/2}\over 1-Q^{j/2}}\times\cr
&(1-e^{2\pi i k}Q^{j})
      (1-e^{-2\pi i k}Q^{j})(1-e^{4\pi i k}Q^{2j-1})
      (1-e^{-4\pi i k}Q^{2j-1})&\eqnu\cr
\label\phiqq\eqnum*
&=\ -{\sqrt{\pi a}\over 2}q^{3k^2\over 2} 
(q^{k\over 2}-q^{-{k\over 2}})
\prod_{j=1}^\infty {1-q^{j/2}\over (1-q^{j})^2}\times\cr
& \prod_{j=1}^\infty
(1-q^{j-k})(1-q^{j+k})(1+q^{j/2-1/4+k/2})(1+q^{j/2-1/4-k/2}).
&\eqnu
\label\phiq\eqnum*
}
$$
As $a\to \infty$, 
$$
\eqalignno{
&\phi(k)\to \sin(\pi k),\ \
(a/4)^{k-1/2}\GG_q^{\sharp}(k)\to \sin(\pi k)\Ga(k+1/2).
&\eqnu
\label\philim\eqnum*
}
$$
\label\GG\theoremnum*
\endproclaim

Concerning the limits, the convergence of $\phi(k)$ to  
$\sin(\pi k)$ is very fast. They practically coincide
even for $a\sim 1$ because $Q$ is getting very close to zero.
The $Q^{1/2}$ is dual to $q^{1/4}=\exp(-1/(4a))$ in the sense of the 
functional equation for the theta function:
$$
\eqalignno{
&\sum_{j=-\infty}^{\infty} Q^{j^2/2}\ =\ (\phi(k)/\sin(\pi k))(k=0)\ =\
{1\over 2\sqrt{\pi a}}\sum_{j=-\infty}^{\infty} q^{j^2/4}.
&\eqnu
\label\phio\eqnum*
}
$$
Here we  first check (\ref\phiqq) and then
use the functional equation 
to go to (\ref\phiq). The latter is straightforward.
Combining (\ref\ggdef) with (\ref\phiq) we can eliminate $Q$.
Indeed, dividing by the  common factors
and distributing $\ \exp({k(k+1)\over 2a})\ $ properly
one  arrives at the $q$-identity:
$$
\eqalignno{
& \sum_{j=0}^\infty q^{j^2-2kj\over 4}{ 1-q^{j+k} \over 1-q^{k}  } 
\prod_{l=1}^j
{1-q^{l+2k-1}\over
 1-q^{l}}\ =&\eqnu\cr
\label\phiden\eqnum*
&\prod_{j=1}^\infty
{(1-q^{j/2})(1-q^{j+k})(1+q^{j/2-1/4+k/2})(1+q^{j/2-1/4-k/2})\over
(1-q^j)}.
}
$$

It is nothing else but the formula (1.11) from [C1] for
the root system $A_1$ and the polynomial spherical representation.
The simplest particular cases  $k=0,1,2,1/2$  readily result
from the  product
formula  for the classical $\vph=\sum_{n=-\infty}^\infty q^{n^2}$
(see (\ref\phio)).
It is possible that
this identity for $k\in \Z_+/2$ 
follows from known formulas.
Then one can extend it analytically
to arbitrary $k$. However here we establish (\ref\phiden) via
(\ref\phiqq).

The key observation is that
$$
\eqalignno{
&\phi(k+1)=-\phi(k),\  \phi(k+2p)=q^{3kp+3p^2}\phi(k)\for
p\equal 2\pi a i.
&\eqnu
\label\phishift\eqnum*
}
$$
The second formula results directly from the definition, the
first one  requires the technique of the shift operators which will
be discussed later. These symmetries and the normalization (\ref\phio)
are sufficient to reconstruct $\phi$ uniquely.

\vskip 5pt {\bf A generalization.} 
A natural setup for the theorem (in one 
dimension) is as follows. 
Let $\{p_m(x)\}$ for $m=0,1,2,\ldots$ be the
Rogers-Askey-Ismail polynomials (see [AI]) in terms of $q^{nx}$ for
$n\in \Z.$ They are $x$-even, pairwise orthogonal with respect to
the pairing $\{f,g\}=$ Const Term $(fg\de_k)$ 
(where $\de_k$ is replaced by the
corresponding Laurent series), and normalized by the condition
$$
\eqalignno{
& p^{(k)}_m=q^{mx}+q^{-mx}+\hbox{\ lower\  powers,\ except\ for\ }  
p^{(k)}_0=1.
&\eqnu
\label\progers\eqnum*
}
$$
The multi-dimensional generalization is due
to Macdonald (see e.g. [M2]).
In this definition, $k\neq -1,-2,-3,\ldots,-m+1$ modulo $2\pi a i\Z.$
For instance,
$$
\eqalignno{
&p^{(k)}_1=q^x+q^{-x},\ p^{(k)}_2=q^{2x}+q^{-2x}+{(1-q^k)(1+q)\over 1-q^{k+1}}.
&\eqnu\cr
\label\ptwo\eqnum*
}
$$

Setting
$$
\eqalignno{
\GG_q^{\sharp}(k;m,n)&\equal
{1\over 2}\int_{\infty+\ep i}^{\infty-\ep i} 
p^{(k)}_m(x)p^{(k)}_n(x)e^{-x^2/a} \de_k(x;q)\ dx\ =\cr
&-{a\pi\over 2}\prod_{j=0}^\infty
{(1-e^{-(j+k)/a})(1-e^{-(j-k)/a})\over
 (1-e^{-(j+2k)/a})(1-e^{-(j+1)/a})}\times
&\eqnu\cr
\label\ggmn\eqnum*
\sum_{j=0}^\infty\ p^{(k)}_m({k+j\over 2})&p^{(k)}_n({k+j\over 2}) 
e^{-(k+j)^2\over 4a} \prod_{l=1}^j
{(1-e^{-(l+2k-1)/a})(1-e^{(l+k)/a})\over
 (1-e^{-(l+k-1)/a})(1-e^{l/a})},
}
$$
we deduce from [C1] the identity:
$$
\eqalignno{
&\GG_q^{\sharp}(k;m,n)\ =\ q^{-{m^2+n^2+2k(m+n)\over 4}}p^{(k)}_m({n+k\over 2})
p^{(k)}_n({k\over 2})\GG_q^{\sharp}(k).
&\eqnu
\label\ggmnid\eqnum*
}
$$

The formula also holds for the $q$-spherical functions, generalizations of
$p^{(k)}_m, p^{(k)}_n$ to  arbitrary complex $m,n$, 
and in the multi-dimensional
setting (see ibid. and [AI]).
Moreover, the Gaussian
$\ga=q^{x^2}=e^{-x^2/a}$ can be replaced by  $(h_l)^\star \ga$
for the $q$-Hermite polynomials
$\{h_l(x)\}$ (the Rogers polynomials at $k=\infty$) and the formal conjugation
$(q)^\star=q^{-1}$. This  is possible (in one dimension only)
because $(h_l)^\star \ga$ are eigenfunctions of the shift operator
(see [BI]). 
We can make $l$ complex too switching
to the $q$-Hermite spherical functions.  So a natural setting 
for Theorem \ref\GG is with  three complex parameters
$\{k,l,m\}.$  Hopefully one can add  more
replacing $\de_k$ by an arbitrary truncation of the theta-function
and involving  $BC_1.$ Also (\ref\ggmnid) can be extended
to nonsymmetric Macdonald polynomials.

We note that both (\ref\ggim) and (\ref\phiden) can be deduced
from known identities (Macdonald, Andrews). The general one
involving Macdonald polynomials (or spherical functions) seems
new.

%
%
%
\vfil
\section{Analytic continuation}

In this  section we  construct meromorphic continuations
of the ``plus-counteparts'' of $\ZZ_q^\sharp$ and  $\ZZ_q^{im}$
and give the formulas for the $a$-limits. We mainly use the Cauchy theorem.
There is another approach based on the shift operator technique, which will
be applied in the next section.

Let us start with the 
``sharp'' function:
$$
\eqalignno{
&\ZZ^{\sharp}_{+q}(k|d) \equal 
{1\over 2i}\int_{\infty+\vep i}^{\infty-\vep i } {\de_k(x;q)\ dx
\over \exp(dx^2/a)+1} \for \cr 
&\vep=\ep_+(a,d)=\min\{\sqrt{{\pi a\over 2d}},{\pi a\over2}\},\ \
d>0, \ \Re k>-1/2.
&\eqnu
\label\zzshd\eqnum*
 }
$$
The integration is with respect to the path $\{C^\vep_- - C^\vep_+\}$
from (\ref\srepm). Recall that it goes from $\infty+\vep i$ to
$\infty-\vep i$ through the origin.
The  function  is well-defined and
has no singularities  in the horizontal
strip $K_+^{\sharp}=\{-2\vep<\Im k<+2\vep\}$ for  $\Re k>-1/2.$

\proclaim{Theorem}
The meromorphic continuation of $\ZZ^{\sharp}_{+q}(k|d)$ to all $k$
is given by the formula
$$
\eqalignno{
&\ZZ_{+q}^{\sharp}(k|d)\ =\ 
-{a\pi\over 2}\prod_{j=0}^\infty
{(1-e^{-(j+k)/a})(1-e^{-(j-k)/a})\over
 (1-e^{-(j+2k)/a})(1-e^{-(j+1)/a})}\times\cr
&\sum_{j=0}^\infty {1-e^{-(j+k)/a}\over 1-e^{-k/a}}
\prod_{l=1}^j
{1-e^{-(l+2k-1)/a}\over
 1-e^{-l/a}}{\exp(kj/a)\over \exp({d(k+j)^2\over 4a})+1}.
&\eqnu
\label\zzssdp\eqnum*
}
$$
The set of its  poles  in the strip $K_+^{\sharp}$
 is $\{-1/2-\Z_+\}$. It also has zeros at all integral $k$.
For any $k\in\C$ except for the poles,
$$
\eqalignno{
&\lim_{a\to \infty} ({a\over 4})^{k-1/2}\ZZ_{+q}^{\sharp}(k|d)\ =\ 
\sin(\pi k)\ZZ_+(k|d),\cr
&  \ZZ_+(k|d)\equal d^{-1/2-k} \ZZ(k)
(1-2^{1/2-k}).
&\eqnu
\label\limzspd\eqnum*
}
$$
\label\ZZSP\theoremnum*
\endproclaim

Formula (\ref\zzssdp) results from Cauchy's theorem. Thanks to 
the $\exp(dx^2/a)$, the convergence is
for all $k$ (except for the poles). 
The zeros of  $\ZZ_{+q}^{\sharp}(k)$ at  $k\in \Z_+$ come from the
product before the summation. It is less obvious for negative integers
(use Theorem \ref\SHIFTWO below).
Actually it is a general property of
the difference Fourier transform (cf. 
(\ref\ggmnid)).  If $\Re k>-1/2$, 
the limit (\ref\limzspd) is due to   
the integral representation. 
For the other $k$, one can  use the shift operators
(the next section). We note that the calculation
of the limit for $\Re k<-1/2$ is equivalent to that
for the imaginary $q$-zeta.

\vskip 5pt  {\bf Imaginary zeta.}
We define it for $\Re k>0$ as follows (cf. (\ref\zzimpl)):
$$
\eqalignno{
&\ZZ^{im}_{+q}(k|d)=\ZZ_{+q}(k|d) \equal \ (-i)\int_{0}^{i\infty} 
(\exp(-dx^2/a)+1)^{-1} \de_k(x;q) dx.
&\eqnu\cr
\label\zzimd\eqnum*
}
$$
We will drop $\{im\}$ till the end of the section.
Its limit as $a\to \infty$ can be readily calculated:
$$
\eqalignno{
&(a/4)^{k-1/2}\ZZ_{+q}(k)\ \to \  
\ZZ_+(k|d).
&\eqnu\cr
\label\zzlimd\eqnum*
}
$$

The first step is to extend this function to small
negative $\Re k$. Let
$$
\eqalignno{
&\ep\ =\ \min\{\sqrt{{\pi a\over 2d}},{1\over4}\},\ 
-2\ep<\Re k<2\ep.
&\eqnu\cr
\label\twoep\eqnum*
}
$$ 
The following function is analytic for such $k$ except for
the poles at $\{\pi a i +2\pi a i\Z\}$ and coincides
with $\ZZ_{+q}(k|d)$ when $\Re k>0$:
$$
\eqalignno{
&\ZZ^{\flat}_{+q}(k|d) 
\equal {1\over 2i}\int_{\ep-\infty i}^{\ep+ \infty i} 
(\exp(-dx^2/a)+1)^{-1}\de_k(x;q) dx+\zz^\pi_0(k),\cr
&\zz^\pi_n(k|d)\equal\ 
{a\pi\over 2}\prod_{j=0}^\infty
{(1-e^{-(j+k)/a})(1-e^{-(j-k)/a})\over
 (1-e^{-(j+2k)/a})(1-e^{-(j+1)/a})}\times &\eqnu\cr
\label\zzimep\eqnum*
\prod_{l=1}^n &
{(1-e^{-(l+2k-1)/a})(1-e^{(l+k)/a})\over
 (1-e^{-(l+k-1)/a})(1-e^{l/a})}
\sum_{m=-\infty}^{+\infty}{1\over \exp({-d(k+n+2\pi a i m)^2\over 4a})+1}.
}
$$

To see this, we replace $(-i)\int_{0}^{\infty i}$ in (\ref\zzimd) by
 $(-i/2)\int_{-\infty i}^{\infty i}$ and then by   
$(-i/2)\int_{\ep-\infty i}^{\ep+ \infty i} +\zz^\pi_0$ provided that
$0<\Re k<2\ep$. Here we applied  Cauchy's theorem.

\proclaim{Theorem}
Let $\xi_l=((2l+1)\pi a i /d)^{1/2},$ where $l=0,1,,\ldots,$
$$
\eqalignno{
&\s \equal \Pi\cup \La\cup\La_+,\
\Pi=\{-\Z_+/2 +\pi a i\Z\}\setminus \{-\Z_+ +2\pi a i\Z\},\cr
&\La\ =\ \{-2\xi_l -\Z_+ +2\pi a i\Z\}\cup
 \{-2\xi_l^* -\Z_+ +2\pi a i\Z\},\cr
&\La_+\ =\  \{2\xi_l -\Z_+ +2\pi a i\Z\}\cup
 \{2\xi_l^* -\Z_+ +2\pi a i\Z\}.
&\eqnu
\label\sbul\eqnum*
}
$$
By $*$ we mean the complex conjugation.
The series 
$$
\eqalignno{
\ZZ^\pi_{+q}(k|d)\ &=\ \sum_{n=0}^\infty \zz^\pi_n(k|d),\ 
\ZZ^\la_{+q}(k|d)=\sum_{l=0}^\infty \zz^\la_l(k|d),\cr
\zz^\la_l(k|d)\ &=\ {a\pi\over 2d}
\Bigl(\xi_l^{-1}\prod_{j=0}^\infty{
(1-e^{-(j+2\xi_l)/a})(1-e^{-(j-2\xi_l)/a}) \over
(1-e^{-(j+k+2\xi_l)/a})(1-e^{-(j+k-2\xi_l)/a})}+\cr
&(\xi_l^{-1})^*\prod_{j=0}^\infty{
(1-e^{-(j+2\xi_l^*)/a})(1-e^{-(j-2\xi_l^*)/a}) \over
(1-e^{-(j+k+2\xi_l^*)/a})(1-e^{-(j+k-2\xi_l^*)/a})}\Bigr)
&\eqnu
\label\zzbul\eqnum*
}
$$
are absolutely convergent for any fixed $k\not\in \s$ with $\Re k<0.$
The difference 
$$
\eqalignno{
&\ZZ^{neg}_{+q}(k|d)\equal
 \ZZ^\pi_{+q}(k|d) -\ZZ^\la_{+q}(k|d)
&\eqnu
\label\zzveed\eqnum*
}
$$
represents  the meromorphic   continuation of $\ZZ_{+q}(k|d)$
to  all  $\Re k<0$ except for the set of
the poles $\p=\Pi\cup\La.$
\label\ZZVEE\theoremnum*
\endproclaim

It is not difficult to check that the poles of  
$\ZZ^\pi_{+q}$ and $\ZZ^\la_{+q}$ apart from  $\p$
have the same residues and cancel each other in 
(\ref\zzveed). The latter set is discrete in contrast to $\s$ which
is dense everywhere (for generic $a$). We will skip the convergence
estimates because the theorem can be reproved using the shift operators
without termwise consideration.

Actually the self-termination of the poles readily follows from
the integral representation (\ref\zzimep). The latter coincides
with $\ZZ^{neg}_{+q}$ and  has no poles
for $0>\Re k>-2\ep$. If the differences of residues of the poles
 from $\La_+$ are nonzero they can be seen 
in this strip (for generic $a$). This argument and the
theorem can be somewhat generalized as follows.

\proclaim{Proposition} Given $r\in \Z_+$,
the meromorphic continuation of the integral $\ZZ^{im}_{+q}$ 
(see (\ref\zzimd)) from the strip $K_r=\{k, \ -r<\Re k<-r+1\}$
to all negative $\Re k$ is given by the series
$$
\eqalignno{
&\ZZ^{-r}_{+q}(k|d)\equal \sum_{j=0}^\infty\, \hbox{sign}\,
(j-r+0.1)\, \zz^\pi_j(k|d)\ -\ \ZZ^\la_{+q}(k|d). 
&\eqnu
\label\zzveedf\eqnum*
}
$$
Its poles ($\Re k<0, r>0$) belong to $\p_r=\{-r+\La\}\cup \Pi\setminus
\{1/2-r+\pi a i\Z\}.$
\label\ZZF\theoremnum*
\endproclaim

This statement will be involved to analyze $ \ZZ^{neg}_{+q}$ as 
$\Re k\to -\infty$ upon the following renormalization:
$$
\eqalignno{
&\widetilde{\ZZ}_{+q}(k|d)=\be(k)\ZZ_{+q}(k|d),\
\widetilde{\ZZ}^{-r}_{+q}(k|d)=\be(k)\ZZ^{-r}_{+q}(k|d)\for 
&\eqnu\cr
\label\zzstd\eqnum* 
&\be(k)\equal
\prod_{j=0}^\infty {(1-q^{j+2k})(1-q^{j+1})\over(1-q^{j+k})(1-q^{j+k+1})}.
&\eqnu
\label\zzbeta\eqnum* 
}
$$
These functions also have  meromorphic continuations to
$k\in \C$. The additional factor 
simplifies the  singularties. The poles of the continuation of 
$\widetilde{\ZZ}_{+q}$ (all $k$) form the set
$\widetilde{\p}=\La\cup\{-1-\Z_+ +2\pi a i\Z\}$. The poles of
the second function (negative $\Re k, \ r>0$) constitute
$\widetilde{\p_r}=\{-r+\La\}\cup \{-1-\Z_+ +2\pi a i\Z\}$. 
It also  has ``two''   zeros 
$\{1/2-r+\pi a i\Z\}$  modulo $ 2\pi a i\Z$.

Similar methods can be applied to the   $\ZZ_q^{im}$ from
Section 1 corresponding to the standard kernel
$(\exp(-x^2/a)-1)^{-1}$.
One has to replace  $\zz^\pi_n(k|d)$
switching to the factors 
$(\exp(-(k+n+2\pi a i m)^2/ (4a))-1)$ and
use $\xi_l=(l\pi a i/d)^{1/2}$ starting the summation from $l=1.$ 
Other ingredients remain the same.

\vskip 5pt {\bf $\bold{k}$-limits.} 
The behavior of $q$-zeta functions as $k\to \infty$
is very non-classical. Let us first tend $\Re k$ to $+\infty$.
We may use  the integral representation:
$$
\eqalignno{
&\ZZ_{+q}(k|d) \to \ \psi_+=(-i)\int_{0}^{\infty i} 
{\prod_{j=0}^\infty (1-q^{j+2x})(1-q^{j-2x})\, dx\over
\exp(-dx^2/a)+1} 
&\eqnu
\label\zzlimplus\eqnum*
}
$$
 as $\Re k\to +\infty$.
So the limit is a positive constant.  

The calculation is more involved when $\Re k$
goes to $-\infty$. Namely, the integral representation cannot be used
anymore and the renormalization (\ref\zzstd) is necessary.
We will apply Proposition \ref\ZZF
as $r\to \infty$ and $k=\ka-r+\nu i$ for $0<\ka<1$. First of all,
$$
\eqalignno{
&
\phi(\nu;\ka) \gets \widetilde{\ZZ}^{-r}_{+q}(k|d)
\to -\widetilde{\ZZ}^\pi_{+q}(k|d) -\widetilde{\ZZ}^\la_{+q}(k|d)
&\eqnu
\label\phidf\eqnum*
}
$$
for a continuous  function
$\phi(\nu;\ka)$ which is $2\pi a$-periodic with respect to
$\nu\in \R$. Here we use the integral representation for
$\widetilde{\ZZ}^{-r}_{+q}(k|d)$ and evaluate the change
of the integrand as  $k\mapsto k-1$
(the multiplicator) for big $-\Re k$.
Thanks to the renormalization, the integral 
tends to a $\Z$-periodic function
of $k.$ 
Of course it is necessary to avoid the
singularities $k\in \s$. 

Thus we established  a relation between 
$\ZZ^\pi_{+q}$ and $\ZZ^\la_{+q}$ in the considered limit. Let us use it:
$$
\eqalignno{
&\lim\widetilde{\ZZ}^{neg}_{+q}\ =\ \lim(\widetilde{\ZZ}^\pi_{+q}-
\widetilde{\ZZ}^\la_{+q})\ =\
-\phi(\nu;\ka)-2\lim\widetilde{\ZZ}^\la_{+q}(k|d).
&\eqnu
\label\phired\eqnum*
}
$$ 
If $k$ is fixed (and is in a general position),
$$
\eqalignno{
&\ZZ^\la_{+q}(k|d) = {1\over 2i}\int_ {\ll}
{\be(k)\de_k(x;q)\, dx\over
\exp(-dx^2/a)+1},
&\eqnu
\label\zzll\eqnum*
}
$$
where the integration contour $\ll$ is the boundary of
a neighborhood of the union  
$$\{(\ep+\R_+)(1+i)\cup (\ep+\R_+)(1-i)\}.$$
The orientation is standard. The  neighborhood must contain no 
points from $\pm \{(k/2+\Z_+/2)+\pi a i \Z \}$. So it depends on $\ka,\nu$
and approaches  the diagonals for big $\Re x$. Now we can go from $k$ to 
$k-1$ evaluating  
$\be(k-1)\de_{k-1}(\be(k)\de_{k})^{-1}$ as $r\to\infty.$

\proclaim{Theorem}
Given $0<\ka<1$ such that $\{\ka+i\R+\Z\}\cap\widetilde{\p}=\emptyset$,
the limit
$$
\eqalignno{
& \psi(\nu;\ka)\ =\ 
\lim_{\Z\ni r\to \infty}\widetilde{\ZZ}^{neg}_{+q}(\ka-r+i\nu)
&\eqnu
\label\psidf\eqnum* 
}
$$
exists and is a continuous $2\pi a$-periodic function of $\nu\in \R.$ 
\label\ZZKLIM\theoremnum*
\endproclaim

The theorem is the  ``plus-variant'' of Theorem \ref\ZZIM, where
$\widetilde{\p}$ was somewhat different.
The set of admissible $\ka$  is dense everywhere in $[0,1]$
because  so is $\widetilde{\p}+i\R+\Z\subset\C$
(for generic $a$). However the function $\psi(\nu;\ka)$ is
well-defined. It is directly  related to the existence of
the totally singular  zeta $\tze_q$ from Section 1,
because $\ZZ^\la_{+q}$ is nothing else but a variant of the function
from (\ref\tzere).

Thanks to the theorem, the asymptotic
number of zeros $N_r$ of  $\ZZ^{neg}_{+q}(k|d)$ for  big (generic) $r$ in
the box
$$B_r\equal \{\ka-r-1<\Re k<\ka-r\} \hbox{\ modulo\ }  2\pi a i\Z$$
coincides with that for  $\widetilde{\ZZ}^{neg}_{+q}(k|d)$ and 
 equals  the number of the poles, i.e.
the points from  $\widetilde{\p}\cap B_r.$ Indeed, it
suffices to calculate the 
difference of the changes
of $\arg(\psi(\nu))$ over the ``periods'' at $\ka=0$ and  $\ka=1.$  
Note that  $N_r$
is quadratic in $r:$ 
$$
\eqalignno{
&N_r= 2\Bigl[{1+d(r+1-\ka)^2\over 4a\pi}\Bigr]+1 , \ 0<\ka<1,\ r\gg 1.
&\eqnu
\label\nfstab\eqnum* 
}
$$

It is likely that 
the positions of individual  zeros (modulo $\Z+2\pi a i\Z$)
approach certain limits  when $r\to\infty$. To ensure this,
the convergence to $\psi(\nu;\ka)$ must be fast enough.
This can be seen
numerically. For instance, if $a=1=d$ the function
$\ZZ_{+q}(k|d)$ has the following (unique?) stable real zero: 
$$
\eqalignno{
& z_6= -6.668819069126,\  z_7= -7.669834033011,\cr
& \hbox{stable\ zero\ } z_r= -r-0.67042591089415.
&\eqnu
\label\stzero\eqnum* 
}
$$
The stabilization is rapid (for small $a$).
At least one  real zero must exist in $B_r$ for sufficiently
big $r$ because $N_r$ is
odd. Recall that  $\ZZ_{+q}(k|d)$ is  real on $\R$, so
the complex zeros appear in pairs.

\vskip 5pt {\bf $\bold{a}$-limits.}
Let us establish the relations to the classical
counterparts of the considered zeta functions:
$$
\eqalignno{
&\ZZ_+(k|d)\ =\  d^{-1/2-k}(1-2^{1/2-k}) \ZZ(k),\cr 
&\widetilde{\ZZ}_+(k|d)\ =\  d^{-1/2-k}(1-2^{1/2-k}) 
\sqrt{\pi}\,\Ga(k+1)\ze(k+1/2).
&\eqnu
\label\zzpnoq\eqnum*
}
$$
See (\ref\limzspd) and (\ref\zzstabl).

\proclaim{Theorem}
Let $\ZZ^{an}_{+q},\ \widetilde{\ZZ}^{an}_{+q}$ be the meromorphic 
continuations of  $\ZZ_{+q}$ and $\widetilde{\ZZ}_{+q}$ to
all $k$ constructed above.
Then
$$
\eqalignno{
&\lim_{a\to\infty}(a/4)^{k-1/2}\ZZ_{+q}(k|d)\ =\ \ZZ_+(k|d) \for
k\neq -1/2,-3/2,\ldots,
&\eqnu\cr
\label\zzpiml\eqnum* 
&\lim_{a\to\infty} a^{k-1/2}\widetilde{\ZZ}_{+q}(k|d) = 
\widetilde{\ZZ}_+(k|d) \for k\neq -1,-2,\ldots.
&\eqnu
\label\zzpstabl\eqnum*
}
$$
\label\ZZPLIM\theoremnum*
\endproclaim

The proof is based on the integral representation (\ref\zzimd)
for $\Re k>0$. As to $\widetilde{\ZZ}_{+q}$, we also use
(\ref\galim) from Lemma \ref\MOAK. When $\Re k>-1/2,$ we can go
 from $\de_k$ to $\mu_k$ 
(see Section 2 and below) or involve $\ZZ^{\flat}_{+q}$ from (\ref\zzimep).
Let us deduce the formula for the limit of $\ZZ^{an}_{+q}$ 
for  $\Re k<-1/2$ from
(\ref\limzspd), Theorem \ref\ZZSP.
Generally speaking,  the real and imaginary 
$q$-zeta functions are not
connected, because the contours of integration are different and
the signs of the Gaussians are opposite.
However, there is a  relation  in the limit.

First, let us consider  $\ZZ^\pi_{+q}(k|d)$ (see (\ref\zzveed)
and (\ref\zzbul)). Since $a\to \infty,$
the terms in the sums (\ref\zzimep) for $\zz^\pi_n(k|d)$ 
with $m>0$  are very small and  can be ignored.
Therefore
$$
\eqalignno{
&\ZZ^{neg}_{+q}+\ZZ^\la_{+q}\ \cong\ 
\ \sum_{j=0}^\infty {a\pi\over 2}\prod_{j=0}^\infty
{(1-e^{-(j+k)/a})(1-e^{-(j-k)/a})\over
 (1-e^{-(j+2k)/a})(1-e^{-(j+1)/a})}\times\cr
&\prod_{l=1}^j
{(1-e^{-(l+2k-1)/a})(1-e^{(l+k)/a})\over
 (1-e^{-(l+k-1)/a})(1-e^{l/a})}
{1\over \exp({-d(k+j)^2\over 4a})+1}\ =\cr
& -\CC_q(k)+\ZZ^{\sharp}_{+q}.
&\eqnu
\label\zzapprox\eqnum*
}
$$ 
Here we used the formula for the ``sharp constant term'' from
(\ref\ccdef) and
of course the identity $(\exp(x)+1)^{-1}+(\exp(-x)+1)^{-1}=1.$

Since $\CC_q^{\sharp}(k)=O(a)$ (Theorem \ref\CT) and
$\ZZ^{\sharp}_{+q}=O(a^{1/2-k}),$ the latter will dominate for big $a$
because $\Re k<-1/2$. Thus 
$\ZZ^{neg}_{+q}\cong \ZZ^{\sharp}_{+q}-\ZZ^\la_{+q},$ and thanks
to (\ref\limzspd) it suffices  to calculate the
limit of $\ZZ^\la_{+q}.$
This can be done via the Moak theorem (Lemma \ref\MOAK)
exactly as we did  in (\ref\limd) with $\tze_q:$
$$
\eqalignno{
&\ZZ^\la_{+q}(k|d)\cong \sum_{l=0}^\infty {a\pi\over 2d}(a/4)^{-k}
((ax_l)^{-1/2} (-x_l)^{k}+(-ax_l)^{-1/2} (x_l)^{k})
&\eqnu\cr
\label\zzbul\eqnum*
}
$$
for $x_l\ =\ ((2l+1)\pi  i /d).$ Then we use the classical
functional equation.

Summarizing,  our technique of
analytic continuation of the imaginary $q$-zeta 
is a $q$-variant of the so-called  Riemann's first proof 
of the functional equation (see [E], 1.6 and the appendix, the
translation of Riemann's ``On the number of primes less than a given
magnitude''). In contrast to the classical
theory, the resulting functions are periodic in the imaginary direction and
have limits as $\Re k\to \pm\infty.$
The ``sharp'' $q$-zeta functions $\ZZ_q^\sharp$ and  $\ZZ_{+q}^\sharp$ 
can be represented
as convergent series for all $k$, which has no counterpart for 
Riemann's zeta.

%
%
%
\vfil
\section{Shift operator}

The shift operator in the setup of this paper
is the so-called Askey-Wilson operator,
the following  $q$-deformation of the differentiation:
$$
\eqalignno{
&Sf(x)\equal (q^x-q^{-x})^{-1}(\tau-\tau^{-1})f(x),\cr 
&\tau f(x)=f(x-{1\over2}),\ q=\exp(-1/a).
&\eqnu
\label\shift\eqnum*
}
$$
For instance,
$$
\eqalignno{
&S\ga \ =\ -q^{1/4}\ga,\ \ S\ga^{-1}\ =\ q^{-1/4}\ga^{-1} \for
\ga=q^{x^2}.
&\eqnu
\label\shiga\eqnum*
}
$$
The name is ``shift operator''  because its action on the hypergeometric
functions, both classical and basic (difference), results in the
shift of the parameters [AW]. For example,
$$
\eqalignno{
&Sp^{(k)}_m(x;k) \ =\ (q^{-m/2}-q^{m/2}) p^{(k)}_{m-1}(x;k+1),\ m=1,2,\cdots
&\eqnu
\label\ship\eqnum*
}
$$
for the Rogers polynomials from (\ref\progers).

The differential shift operators for arbitrary
root systems are due to Opdam (see [O] and also [H] for
the interpretation via the Dunkl operators). The general
difference ones are considered  in [C3]. They 
depend on $k$ in contrast to the simplest case considered
here ($A_1$). The celebrated constant term conjecture
and the Macdonald-Mehta conjecture were verified using these operators
(see [M1, O, C3, C1]). They were also used in [O] for analytic
continuations. 

In this section, the shift operators will be applied to go from 
positive $\Re k$ to the left without involving Cauchy's theorem
and convergence estimates of the pole decomposition.
Some properties of the difference shift operators have no counterparts in
the classical theory. For instance, they give a connection between
the values of $q$-zeta functions at $k$ and $k+1$ which collapses in the
limit. We do not discuss the technique of the shift opeartors
in detail.  It is similar to  the algebraic 
theory [C1,C3]. See also [C2]. 

\vskip 5pt {\bf Shift-formula.}
Let $\pi a>\ep>0$ and  $\int$ be one of the following integrations:
$$
\eqalignno{
&\int_{ct} = {1\over \pi a i}\int_{-\pi a i/2}^{\pi a i/2}dx,\ \
\int_{im} =  {1\over 2i}\int_{\ep-\infty i}^{\ep+\infty i}dx,\cr
&\int_{re} = 
{1\over 2}\int_{-\infty -\ep i}^{+\infty -\ep i}dx,\ \ 
 \int_{\sharp} = {1\over 2i}\int_{\infty+\ep i}^{\infty-\ep i }dx.
&\eqnu
\label\ints\eqnum*
}
$$

Given $u>0$,  a function  $E(x)$ is called $u$-regular,
if it is even,  analytic in
$$
\eqalignno{
&D_{ct} = \{x|-u\le \Re x\le u\} \for \hbox{CT}, \cr
&D_{im} = \{x|-u-\ep\le \Re x\le u+\ep\} \for \int_{im},
&\eqnu\cr
\label\doms\eqnum*
&D_{re} = \{x\in \pm\ep i+\R\} \for \int_{re}\ \ (\hbox{no\ } u),\cr 
&D_{\sharp} = D_{re}\cup
\{x|-\ep\le \Im x\le \ep,\ -u\le \Re x\le u\} \for \int_{\sharp},
}
$$
and continuous on the  closure of $D.$
It is also assumed to be  $\pi a i$-periodic for CT and  have
coinciding continuous  limits
$$e^{\pm}(\ka+\nu i)\ =\ 
\lim_{\nu\to \pm \infty} E(x)\de_{\ka+\nu i}(x)
$$
in the case of the imaginary integration $\int_{im}$. 

The following proposition results  from (\ref\ship).
The proof  is similar to   [C2] (in the one-dimensional
setup).

\proclaim{Proposition} Let $K_\ep$ be the set of $k$ such that
  $\Re k>2\ep$ for the imaginary
integration, $\Re k>0$ for CT,  $\{\Re k>0,\ |\Im k|<2\ep\}$ 
for the sharp one
$\int_\sharp$, and  $K_\ep= \{\Im k+2\pi a\Z\neq \pm 2\ep\}$ 
for the real integration.   Given a
$1/2$-regular $E$ and  
a Rogers polynomial $p_m^{(k)} (m>1),$
$$
\eqalignno{
&\int S(E) p^{(k+1)}_{m-1}\de_{k+1}\ =\ 
q^{k}(q^{-k-{m\over2}}-q^{k+{m\over 2}})\int E p^{(k)}_{m}\de_{k}
&\eqnu
\label\shiftp\eqnum*
}
$$
for $k\in K_\ep,$ provided the existence of the integrals.
\label\SHIFTP\theoremnum*
\endproclaim

Note that we do not take 
$i\R$ as the path for the imaginary
integration because  $S(E)$ may have singularities
at $2\pi a i\Z.$
Concerning the convergence of the integrals from (\ref\shiftp),
the function $E$ can be arbitrary (say, continuous) for CT. For
the imaginary integration, we can use the Stirling formula
from Lemma \ref\MOAK. 
For instance, 
the integrability  of  $|E(x)||x|^{-2\Re k}$ in the imaginary
directions inside $D_{im}$ is sufficient.
The estimates in the real direction are more involved but
not difficult either. 

Next, we  put 
$$
\eqalignno{
&\de_{k+1}\ =\ (1-q^{2x+k})(1-q^{-2x+k})\de_k,
&\eqnu
\label\dekpl\eqnum*
}
$$
express $ (1-q^{2x+k})(1-q^{-2x+k})$ in terms of $p^{(k)}_2$ 
(see (\ref\ptwo)), and then apply 
the proposition twice to turn  $p^{(k)}_2$
into a constant.

\proclaim{Theorem} Given  $\ep>0,$ 
let $E(x)$ be $1$-regular. Assuming (again) that $k\in K_\ep$
and the integrals are  well-defined,
$$
\eqalignno{
&{(1-q^{2k+1})(1-q^{2k})\over (1-q^{k+1})(1-q^{k})}
\int E\de_{k}=\cr
&\int E\de_{k+1}\ +\ {q^{3/2+k}\over (1-q^{2k+2})(1-q^{2k+3})}
\int S^2(E)\de_{k+2}.
&\eqnu
 \label\shiftwo\eqnum*
}
$$
\label\SHIFTWO\theoremnum*
\endproclaim

The generalization of the theorem to arbitrary root systems is
straightforward. However the explicit formulas for $\de_{k+1}/\de_k$ in
terms of ``small'' Macdonald's polynomials are getting  
complicated. 

An immediate application of the theorem is
the constant term conjecture in the $A_1$-case (see [AI,M1]
and the papers of Milne and Stanton from [A]).
Let $\int=\int_{ct},\ E=1.$ Then
$$
\eqalignno{
&\int_{ct}\de_{k}\ =\ {(1-q^{k})(1-q^{k+1})\over (1-q^{2k})(1-q^{2k+1})}
\int_{ct}\de_{k+1}=\ \cdots\ =\cr
&2\prod_{j=0}^\infty {(1-q^{k+j})(1-q^{k+j+1})\over (1-q^{2k+j})(1-q^{j+1})}.
&\eqnu
\label\shifct\eqnum*
}
$$

The resulting product formula holds for $\Re k>0.$ 
The reasoning is the same for other integrations with a reservation
about the convergence conditions.
For instance, we 
get (\ref\ccprod) for the sharp integration. However now
the convergence is for $\Re k<0,$ so we need to employ a variant of the
theorem for the negative countours of integration and use
the chain
$$
\eqalignno{
&\int \de_{k}\ =\ {(1-q^{2k-1})(1-q^{2k-2})\over (1-q^{k})(1-q^{k-1})}
\int \de_{k-1}=\ \cdots\ .
&\eqnu
\label\shifneg\eqnum*
}
$$
Actually this chain can be applied to establish (\ref\shifct)
as well. First, we assume that $k\in \Z_+$, second, get the exact
formula and, third, use the
analytic continuation.

It is not surprising that the products in
 (\ref\shifct) and  (\ref\ccprod) coincide.
Formally, 
$$-\pi a i\int_{ct}E \de_{k}\ =\  2i\int_{\sharp}E \de_{k}
$$
if $E$ is regular everywhere.

Another application is to the Gaussians. Using (\ref\shiga)
we come to (\ref\ggim) for $E=\ga^{-1}.$ When the integration
is sharp and  $E=\ga,$ we establish
the symmetry (\ref\phishift), which played the key role 
in Theorem \ref\GG. 

Let us apply $S^2$ to our real and imaginary
plus-kernels $E^{(\pm)}=(1+q^{\pm x^2})^{-1}$
when $d=1:$
$$
\eqalignno{
&S^2(E^{(\pm)}) = \pm
{q^{x^2+1/2}(q^{x^2+1}-1)\over 
(1+q^{x^2})(1+q^{(x+1)^2})(1+q^{(x-1)^2})}.
&\eqnu
\label\shiftpm\eqnum*
}
$$
So it is the same for the imaginary and real integration up to a sign.
The complete factorization is a special feature of $S(E^{(\pm)})$
and $S^2(E^{(\pm)}).$ The formulas are  not  that nice for
$d\neq 1$ and for higher $S$-powers of $E^{(\pm)}$. We mention that 
big $S$-powers  have remarkable stabilization
properties which we hope to discuss somewhere else.

\vskip 5pt {\bf Analytic continuation.}
The theorem can be applied infinitely many times if
the function $E$ is $\infty$-regular, for instance, 
in the case of the real integration $\int_{re}$.
Of course it is necessary to ensure the convergence of the 
involved integrals
and the resulting series. It may lead to interesting transformations
and identities. However it seems that the most promising  applications
are analytic. Namely, using (\ref\shiftwo), we can extend the function
$\int E \de_k$ to negative $\Re k$ when
$E$ is analytic in the entire complex
plane or proper vertical strips.
If $E$ is meromorphic the method can still be used.
Of course the poles will contribute to (\ref\shiftwo), but the extra
terms are meromorphic  for reasonable 
$E.$

\proclaim{Corollary}
In the setup of Theorem \ref\SHIFTWO,
let $\EE(k)\equal\int E(x)\de_k$ for $k\in K_\ep.$
The function
$$
\eqalignno{\
\EE^{(1)}(k)\ =\ 
{(1-q^{k+1})(1-q^{k})\over
(1-q^{2k+1})(1-q^{2k})}
\EE(k+1)\ +&
&\eqnu\cr
\label\eeone\eqnum*
{q^{3/2+k}
(1-q^{k+2})(1-q^{k+1})\over
 (1-q^{2k+4})(1-q^{2k+3})(1-q^{2k+2})(1-q^{2k+1})}
&\int S^2(E)\de_{k+2}
}
$$
is a meromorphic  continuation of $\EE(k)$ to $K_\ep -1,$
with $\int S^2(E)\de_{k+2}$ being analytic.

Assuming that $E$ is $r$-regular for $r\in 1+\Z_+$ and
all the integrals are convergent,
the iterations of (\ref\eeone) provide a meromorphic
continuation of $\EE(k)$ to $K_\ep -r$ with simple
poles belonging to the set 
$$\{-\Z_+/2+\pi a i\Z\}\setminus\{-\Z_+ +2\pi a i\Z\}.$$
\label\EE\theoremnum*
\endproclaim

 The described procedure becomes somewhat sharper for the
imaginary integration when $\de_k$ is replaced by
$\mu_k$ from (\ref\mu). Namely,
$$
\eqalignno{\
&\EE(k)\ = \ 
 2(1+q^k)^{-1}\int_{im} E(x)\mu_k
&\eqnu
 \label\eeint\eqnum*
}
$$
for $k\in K_\ep,$
where the latter integral gives a meromorphic continuation of $\EE(k)$ 
to $K_\ep -1/2$ with the poles in $\pi a i+2\pi a i \Z.$
Respectively, substituting $\mu$ for $\de$, the corollary
provides the continuation to  $K_\ep -r-1/2$ for
$r$-regular $E$. 

Note that the above  set of poles coincides with $\Pi$
from (\ref\sbul). Thus the results of the 
previous section can be reproved by
means of the  shift operator technique  for
 $E(x)=$ $(\exp(-dx^2/a)+1)^{-1}$ with a 
reservation about meromorphic $E.$
Using (\ref\eeone) only, we have to stop
once the first pole of $E$ appears, i.e. when  
 $\Re k=-2\Re(\xi_0)$
(Theorem \ref\ZZVEE).

\vskip 5pt   {\bf Limits.} An important application of the 
shift operator is to 
the $k$-limits and $a$-limits. We will discuss here only the
latter. Formula (\ref\limzspd) is still
not justified for  $k\le -1/2$.
Now we can proceed as follows. Let us take $E=$ $(\exp(dx^2/a)+1)^{-1}.$

The diagonal poles are departing from the origin as $a\to \infty,$
so we can define the analytic continuation
using $\EE, \EE^{(1)}, \EE^{(2)}$ etc. (any number
of them). It suffices to consider the first
step, the passage from $\Re k>-1/2$ to  $\Re k>-3/2.$
Then we can repeat the procedure.
For big $a,$
$$
\eqalignno{\
& ({a\over 4})^{k-1/2}\EE^{(1)}(k)\ \cong\ 
{(k+1)\over a(k+1/2)}({a\over 4})^{k+1/2}
\EE(k+1)\ +\cr
& {1\over
 (k+3/2)(k+1/2)}
 ({a\over 4})^{k+3/2}\int_{\sharp} S^2(E)\de_{k+2}.
&\eqnu
 \label\eehalf\eqnum*
}
$$
\comment
$$
\eqalignno{\
& ({a\over 4})^{k-1/2}\EE^{(1)}(k)\ \cong\ 
({a\over 4})^{k-1/2}
{(k+1)\over 4(k+1/2)}
\EE^{(1/2)}(k+1)\ +\cr
& {2  ({a\over 4})^{k-1/2}
{a^2\over
 8(2k+3)(2k+1)}
\int_\sharp S^2(E)\mu_{k+2},\  
\int_{\sharp}\ =\ {1\over 2i}\int_{\infty+\ep i}^{\infty-\ep i }dx
&\eqnu
 \label\eeone\eqnum*
}
$$
\endcomment

We already know (thanks to the initial integral
representation) that 
$$
({a\over 4})^{k+1/2}\EE(k+1)\to \sin(\pi k)\ZZ_+(k+1|d)
\for k>-3/2. 
$$
Actually the exact formula does not matter since this term comes
with  the factor $(1/a)$ and will vanish.
So we need to evaluate the second term only. For the sake of
simplicity, let $d=1.$ Using  (\ref\shiftpm),
$$
\eqalignno{ 
&\lim_{a\to\infty}({a\over 4})^{k-1/2}\EE^{(1)}(k)\ =\cr 
&{1\over 2i(k+3/2)(k+1/2)}\int_{\infty+\vep i}^{\infty-\vep i } 
z^{-1/2}(-z)^{k+2}\e^{(2)} dz,\cr
&\e^{(2)}\ =\ 
-{e^{z}(1-e^{z})\over 
(1+e^{z})^3}.
&\eqnu
 \label\eelimtwo\eqnum*
}
$$

Here we substituted $z=x^2/a.$ The integration path must be of course
recalculated. However it  can be deformed to
the same ``sharp'' shape for a small $\vep>0$ because  $k+2>1/2.$ 
In fact, $k+2>-1/2$ is sufficient. Then we may integrate by parts:
$\e^{(2)}$ is the second derivative of $\e=(1+e^z)^{-1}.$
The latter is not surprising because $S^2$ is a deformation
of the second derivative. Finally,
$$
\eqalignno{ 
&\lim_{a\to \infty}
({a\over 4})^{k-1/2}\EE^{(1)}(k)\ =\  
{1\over 2i}\int_{\infty+\vep i}^{\infty-\vep i } 
z^{-1/2}(-z)^{k}\e dz\ =\cr
&(1-2^{1/2-k})\sin(\pi k) \ze(k+1/2)\Ga(k+1/2)=\sin(\pi k)\ZZ_+(k|1).
&\eqnu
 \label\eelim\eqnum*
}
$$

The last equality holds for $\Re k>-1/2.$ The integral 
from (\ref\eelimtwo) represents $\ZZ_+(k|1)$ for  $\Re k>-5/2$
(we need it only when  $\Re k>-3/2$). So, indeed,
we moved by $-1$ to the left. Similarly, we can use
$\EE^{(2)}$ and so on. This proves   
(\ref\limzspd). One  can  (try to)
apply this method for finding other terms of
the $(1/a)$-expansion of  $\ZZ^\sharp_{+q}(k|d)$
and other $q$-zeta functions. Moreover, 
it may  work as $a\to 0$ and $q$ tends to roots of unity,
however this is beyond the framework of the paper.

Summarizing, the technique of analytic continuation from this section
has the following  classical origin. 
Each integration by
parts shifts by $-1$ the $k$-domain of the 
integral representation (\ref\eelim) for the ``plus-zeta''.
There are of course other similar examples
in the classical theory of zeta. However the
difference integrating by parts, namely (\ref\eeone), seems more promising.
It has better analytic properties and, what is important,
establishes a certain connection between the values of $q$-zeta
functions (or more general integrals) at $k$ and $k+1,$
collapsing in the limit. The stablity of our zeta functions
modulo $\Z,$ which has no classical counterpart, is also directly
related to  (\ref\eeone).

%
%
%
\vfil
\section{Locating the zeros}
In this section we discuss (mainly numerically)
the zeros of $q$-zeta functions and demonstrate that their 
behavior is quite
regular at least for $\Re k>-1/2.$
Almost certainly the evolution of
the $q$-zeros is  related to
the distribution of the classical ones. We will give some evidence.

\vskip 5pt   {\bf Imaginary zeta.} 
Let us begin with the ``plus-imaginary'' $\ZZ^{im}_{+q}(k|d)$
from (\ref\zzimd). Recall that 
$$ 
\eqalignno{
&\ZZ^{im}_{+q}(k|d)\ =\ (-i)\int_{0}^{\infty i}\ 
{1\over e^{-dx^2/a}+1}\de_k(x;q) dx \for\Re k>0, &\eqnu\cr
\label\zzplnum\eqnum* 
\lim_{a\to\infty}&(a/4)^{k-1/2}\ZZ^{im}_{+q}(k|d)=
 \ZZ_+(k|d)=(1-2^{1/2-k})d^{-1/2-k} \ZZ(k),\cr
& \hbox{where\ } \ZZ(k)\ =\
\ze(k+1/2)\Ga(k+1/2).
}
$$
Its analytic continuation has the same limit for all $k$
except for the poles. Since $\ZZ^{im}_{+q}$ is 
$2\pi a i$-periodic and has the limit (a positive constant)
as $\Re k\to\infty$ the number of its zeros $N_+^\ep(a;d)$
in the strip $K_+^\ep=\{k\ |\ \Re k >\ep>0\}$ 
modulo $2\pi a i$ is as follows:
$$ 
\eqalignno{
&N_+^\ep(a;d)\ =\ 
{1\over \pi}\arg(\ZZ^{im}_{+q})\, \vert_{k=\ep}^{\ep+\pi a i}.
&\eqnu\cr
\label\zznep\eqnum*
}
$$
Here we used that  $\ZZ^{im}_{+q}$ is real on $\R_+$ and 
$\pi a i+\R_+.$

Numerically, it has  no zeros 
in the right half-plane at least for $a=25,\ \ep=0.05,\ d=2.$
Here is the table of the arguments 
$\al(n)=\arg(\ZZ^{im}_{+q})(\pi a i n/21+\ep)$ for $n=1,\ldots,
20$ ($\al(0)=0=\al(21)$):
$$
\eqalignno{
\al(1-5)&=-7.45845\  -11.6586\ -12.1104\ -12.4365\  -11.5873
&\eqnu\cr
\label\altable\eqnum*
\al(\ 6-10)&=-10.1478\  -8.26116\ -6.03690\ -3.56672\   -1.01060\cr  
\al(11-15)&=\ \ \ 1.61364\  \ \ \  4.19644\ \ \ \ \, 6.66558\ \ \ \ \ 
8.91015\ \ \ \ \  10.8286\cr
\al(16-20)&=\ \ \ 12.2974\  \ \ \  13.1573\ \ \ \ \, 13.1854\ 
\ \ \ \ 12.0256\  
\ \ \ \ 8.96469.
}
$$
The total number of points on the half-period 
is $\sim 320$ in this calculation.
The behavior of the argument becomes  more complicated 
approximately after $a=30.$ Our computer program (it adds points
automatically to ensure the required accuracy)
was not able to reach a reliable answer for such $a.$
Eventually the zeros will appear in the right half-plane.
Indeed, the $q$-deformations of the
classical zeros with $\Re k=1/2$ coming from the factor $(1-2^{1/2-k})$
must arrive. 

Without going into detail, let us mention that similar results
hold when $(\exp+1)^{-1}$ is replaced by regular modular functions
(cf. (\ref\thimzl)).
The ``imaginary'' $q$-deformations of 
the corresponding $L$-functions (with proper gamma-factors)
are well defined for $\Re k>0.$
Sometimes elementary functions can be used instead of modular
ones in the integral representations. Say, for the $L$-function  
associated with the  ``standard''
Dirichlet character modulo $3$,
we can proceed  as follows:
$$ 
\eqalignno{
&\ZZ^{im}_{3,q}(k|d)\ =\ (-i)\int_{0}^{\infty i}\ 
{1\over e^{-dx^2/a}+1+ e^{dx^2/a}}\de_k(x;q) dx 
\for\Re k>0. &\eqnu\cr
\label\zzdiri\eqnum* 
}
$$
Its analytic continuation 
approaches the classical $L$-function
up to a gamma-factor for all $k\in\C$ as $a\to\infty.$ 
It has no  $q$-zeros in the right half-plane
at least when $a\sim 25, d=2.$

Recall that in the critical strip  $1/2>\Re k>-1/2,$ we use the analytic
continuation, namely, either 
$$
\eqalignno{
& \ZZ^{\flat}_{+q}\ =\ \ZZ^{neg}_{+q} \hbox{\ from\ (\ref\zzimep)}, \or \cr 
&\ZZ^{\mu}_{+q}\ =\
-i(1+q^k)^{-1}\int_{1/4-\infty i}^{1/4+\infty i}\ 
(\exp(-dx^2/a)+1)^{-1}\mu_k(x;q) dx.
&\eqnu
\label\zzannum\eqnum* 
}
$$ 
They coincide (see (\ref\zzveed) and  (\ref\eeint)). The passage to
negative $\Re k$ is  necessary 
because the most interesting zeros, the $q$-deformations
of the classical ones, prefer the  left half-plane.
The theoretical  explanation is as follows. 

Given a zero $k=z$ of the classical $\ze(k+1/2)$ (only
$z\in i\R_+$ will be considered) we can calculate the
linear approximation  $\tz_+(a)$ for the corresponding exact
zero  $z_+(a)$ of  $\ZZ^{\flat}_{+q}:$
$$
\eqalignno{
& \tz_+(a)=z(1+{ 4d^{-1}(z+1/2)\ze_+(z+3/2)-d(z-1)\ze_+(z-1/2)
\over 12 a\ze_+'(z+1/2)}) \cr
&\for \ze_+'(s)=\partial \ze_+(s)/\partial s,\ \ze_+(s)=(1-2^{1-s})\ze(s).
&\eqnu
\label\tzapnum\eqnum* 
}
$$
We  expand the integral representation in terms of $(1/a)$
using a refined variant of Lemma \ref\MOAK. Another approach is
based on the shift operator technique. However we did not check all
the estimates in
(\ref\tzapnum) so the exact range (both $a, k$) where it can
be used  is not determined  at the moment.

The first $\tz_+(a)$ which we found in the right half-plane 
corresponds to $z= 1977.2714i.$ The formula predicts that this $z$
could have the $q$-deformation $z_+(a)$  
to the right, but we did not even try to find it numerically.
It is  well beyond the capacity
of our computer program and the speed of existing (available)
computers. The applicability of the linear formula for so big $z$ 
is also unclear. A more systematic analysis of the Moak-Stirling
formula is necessary.

The negativity of $\Re (\tz_+(a))$ is 
a certain  property of Riemann's zeta and its derivative.
As we will see, it is connected with known facts.
However there is no explanation why the tendency
is that strong.

\vskip 5pt   {\bf Sharp zeta.}
Let us switch entirely to the main object of our computer simulation,
the plus-sharp $q$-zeta function $\ZZ^\sharp_{+q}(k|d)$ from (\ref\zzssdp).
It tends to 
$(a/4)^{1/2-k}\sin(\pi k)\ZZ_+(k|d)$ for all $k$ except for the poles
(Theorem \ref\ZZSP).
We  will  analyze it in the strip  
$$
\eqalignno{
&K_+^{\sharp}\ =\ \{-2\vep<\Im k<+2\vep\} \for 
\vep=\sqrt{{\pi a\over 2d}},
&\eqnu
\label\stripnum\eqnum*
 }
$$
where it has  poles at $\{-1/2-\Z_+\}$ and vanishes at
all integral  $k.$  We assume that $a>2/(\pi d^2)$
which provides that the poles with $\Im k= \pi a$ will
not appear in this strip.

The formula for the linear  approximation  $\tz^\sharp_+(a)$
for the exact zero  $z^\sharp_+(a)$ of  $\ZZ^{\sharp}_{+q}$
corresponding to a given  zero $z$ of $\ze(k+1/2)$ reads
$$
\eqalignno{
& \tz^\sharp_+(a)=z(1-{ 4d^{-1}(z+1/2)\ze_+(z+3/2)-d(z-1)\ze_+(z-1/2)
\over 12 a\ze_+'(z+1/2)}). 
&\eqnu
\label\tzasnum\eqnum* 
}
$$
The reservation about
the applicability is the same as above. 

We will give the numerical values for the zeros $z^\sharp$ deforming
the classical ones $z$
in the strip $K_+^\sharp$ and the values
of the corresponding linear approximations $ \tz_+^\sharp.$
Let $d=2, \ a=750:$ 
\comment
the case of $d=2, \ a=1000:$ 
$$
\eqalignno{
z=&14.1347i\ \ z_+^\sharp=0.0978 + 14.1427i\ \ \
\tz_+^\sharp=0.0977 + 14.1435i\cr
& 21.0220i\ \ \ \ \ \ \  0.2634 + 21.0595i\ \ \ \  \ \ \ 0.2628 + 21.0634i\cr
& 25.0109i\ \ \ \  \ \ \ 0.4249 + 24.9727i\ \ \ \  \ \ \ 0.4309 + 24.9760i\cr
& 30.4249i\ \ \ \  \ \ \  0.681\ + 30.409i\ \ \ \  \ \ \ \ \, 
0.685\ + 30.412i\ \cr
& 32.9351i\ \ \ \  \ \ \  0.829\ + 33.019i\ \ \ \ \  \ \ \ \, 
0.825\ + 33.048i\  \cr
& 37.5862i\ \ \ \  \ \ \  1.248\ + 37.908i\ \ \ \ \  \ \ \ \, 
1.326\ + 38.039i\  \cr
& 40.9187i\ \ \ \  \ \ \  1.44\ \ + 40.84i\ \ \ \ \ \ \,  \ \ \ \, 
1.44\ \ + 40.82i\ \ \cr
& 43.3271i\ \ \ \ \ \ \  1.74\ \ + 43.28i\ \ \ \ \ \  \,  \ \ \ \,
1.84\ \ + 43.32i\ \ \cr
& 48.0052i\ \ \ \ \ \ \  2.2\ \ \ + 47.9i\ \ \ \ \ \ \  \   \ \ \ \  
2.3\ \ \ + 47.7i\ \ \ \cr
& 49.7738i\ \ \ \ \ \ \  2.5\ \ \ + 50.1i\ \ \ \ \ \ \  \   \ \ \ \ 
2.4\ \ \ + 50.5i\ \ \ \cr
& 52.9703i\ \ \ \ \ \ \  2.7\ \ \ + 52.7i\ \ \ \ \ \ \  \   \ \ \ \ 
3.2\ \ \ + 52.1i\ \ \ .
&\eqnu
\label\zertab\eqnum*
}
$$
iterations: first three
\endcomment
$$
\eqalignno{
z=&14.1347i\ \ z_+^\sharp=0.1304 + 14.1450i\ \ \
\tz_+^\sharp=  0.1302 + 14.1465i\cr
& 21.0220i\ \ \ \ \ \ \ 0.3514 + 21.0702i \ \ \ \ \ \ \ 0.3504 + 21.0771i \cr
& 25.0109i\ \ \ \ \ \ \ 0.5641 + 24.9586i \ \ \ \ \ \ \ 0.5745 + 24.9643i \cr
& 30.4249i\ \ \ \ \ \ \ 0.9046 + 30.4014i \ \ \ \ \ \ \ 0.9134 + 30.4077i \cr
& 32.9351i\ \ \ \ \ \ \ 1.1051 + 33.0341i \ \ \ \ \ \ \ 1.0998 + 33.0854i \cr
& 37.5862i\ \ \ \ \ \ \ 1.6449 + 37.9660i \ \ \ \ \ \ \ 1.7675 + 38.1895i \cr
& 40.9187i\ \ \ \ \ \ \ 1.9080 + 40.8119i \ \ \ \ \ \ \ 1.9141 + 40.7816i \cr
& 43.3271i\ \ \ \ \ \ \ 2.2860 + 43.2485i \ \ \ \ \ \ \ 2.4497 + 43.3138i \cr
& 48.0052i\ \ \ \ \ \ \ 2.9259 + 47.8424i \ \ \ \ \ \ \ 3.1103 + 47.5578i 
&\eqnu
\label\zertab\eqnum*
}
$$
We did not try to reach  high accuracy. However the existence of
the  zeros was carefully justified by means of  
sequences of integrations
over  diminishing  rectangles around the consecutive approximations. 
The same integrations were  used to calculate 
the $q$-zeros themselves.

Note that the standard Newton method and its modifications work for the first
5 zeros but diverge for the other 4 even when applied to the 
zeros we found via the contour integrating ($ z_+^\sharp$). 

For $a$ till $10000,$ there is no problem with
finding the first zeros, but calculating the last ones 
is difficult.

The zeros exist for small $a$ as well. If $d=2$ the least $a$
with the zero in the strip $K_+^\sharp$ is $a=61$ (anyway $a=60$ 
is smaller than necessary). The zero 
equals $z_+^\sharp=1.1721 + 13.8088i$ whereas $2\vep=13.8433.$
The corresponding  linear approximation
$\tz^\sharp_+=1.60159 + 14.2789i$ is beyond the $K$-strip.

This does not mean that we cannot consider the  $q$-deformations of the 
classical zeros for small $a.$
For instance, tracing $1.1721 + 13.8088i$
back in $a:$ 
$$
\eqalignno{
a=30:\ &z^\sharp_+\ =\, -3.6951 + 10.5620i,\ 
\tz^\sharp_+= 3.2566 + 14.4279i,\cr
a=1.98952:\  &z_+^\sharp\ =\ -11.5000 + 2.5004i.
&\eqnu
\label\trdtwo\eqnum* 
}
$$
The first is very far from its linear approximation $\tz^\sharp_+,$
the second seems to have nothing to do with the classical
$z=14.1347i$ at all. However both are indeed the
$q$-deformations of the latter.

We note that outside  $K_+^\sharp,$
the appearance of diagonal poles (from $\La$) changes the picture of zeros
(and their number) dramatically.
For instance, $2\vep =9.70813$ for $a=30,d=2$ and 
there is a zero  $0.3862 + 9.9175i.$
It is also necessary to take into consideration the
classical zeros due to the  factor  $(1-2^{1/2-k}).$
For instance, $ 0.8223 + 9.0518i$,\
 $1.9601 + 17.6542i$
are exact deformations of 
$ 0.5 + 9.0647i$,\  $0.5 + 18.1294i$ for $ d=2, a=100$.

According to our calculations, 
it is possible that all zeros $z_+^\sharp$
inside $K_+^\sharp$ (or a bit smaller strip)
and with $\Re k>-1/2\ (a>2/(\pi d^2))$ 
are $q$-deformations of the
classical ones ($\Re k=0, 1/2$). However this must be checked
more systematically. Vice versa, all classical zeros  
presumably  have the  deformations for all $q,$
but  maybe far from the critical strip (see (\ref\zertab), (\ref\trdtwo)).

\vskip 5pt   {\bf Why do the zeros move to the right?}
The last topic we are going to discuss is a qualitative analysis
of the positivity of $\Re (\tz^\sharp).$ A similar question is
about ups and downs of the linear deformations.
They look random in the  table (\ref\zertab), but there are interesting
tendencies, especially for the usual ``minus-zeta''. We will not touch it
upon here.

Since $z\in i\R_+,$
we need to examine the positivity of
$$
\eqalignno{
&\eta_+(z|d)\equal
&\eqnu\cr
\label\etaim\eqnum*
\Im\, \Bigl(\,&{4d^{-1}(z+{1\over 2})(1-2^{-{1/ 2}-z})\ze(z+{3\over  2})
-d(z-1)(1-2^{{3/2}-z})\ze(z-{1\over 2})
\over 12 a\ze'(z+1/2)(1-2^{1/2-z}) }\, \Bigr).
}
$$

If $d=2$ the smallest zero $z$ when (\ref\etaim)
becomes negative is  $z=1977.27i.$
The next one is  $2254.56i.$ 
Actually this property has almost
nothing to do with the zeros of $\ze.$ Indeed, the first
interval where $\eta_+(z|2)$ takes negative values ($\Im z>6$)
is  $[1977.24i, 1977.35i].$ Since there are many zeros of Riemann's
zeta, there is nothing unusual that  one of them appeared in this interval.
The tendency is practcally the same from $d=0.2$ till $d=\infty.$
Only around  $d=0.1,$ a (much) smaller ``negative'' zero appears: 
$163.03i$. The next one is $353.49i.$

Thus the positivity phenomenon is mainly because 
of  the coefficient of $d$ in (\ref\etaim):
$$
\eqalignno{
&\Im\Bigl({
-(z-1)(1-2^{3/2-z})\ze(z-1/2)
\over 12 a\ze'(z+1/2)(1-2^{1/2-z}) }\Bigr).
&\eqnu
\label\etaimp\eqnum* 
}
$$

Here  the factor $(1-2^{3/2-z})/(1-2^{1/2-z})$ 
is not too significant. Numerically, it is because its real part,
always positive, provides the major contribution.
Let us drop it (together with $12a$). So we need to interpret the 
tendency
$$
\eqalignno{
&\eta(z)\equal\Im((1-z)\ze(z-1/2)/ \ze'(z+1/2))>0 
&\eqnu
\label\etaimz\eqnum* 
}
$$
on the zeros $z\in i\R_+$  of $\ze.$

Note that $\eta_+(z|d)$ without the factors $(1-2^{\cdots})$ 
describes the linear deformations in the case of the 
usual zeta (see  (\ref\tzan)). Therefore a little influence of the skipped
factors could be expected a priori.
Actually we  can simlify $\eta$ even more replacing  $(z-1)$ by $i,$
since it almost belongs to $\R_+ i$. However this factor is necessary
to fix the positivity for small zeros, so we will not touch it.

The phenomenon  is again not about the zeros.
The first interval where $\eta(z)$ 
takes negative values is  
$[1267.47i, 1267.70i]$  whereas the first zero with negative $\eta$
is $1267.57i.$  However we will  use that $z$ is a zero
in the following calculation.

Applying the functional equation 
$\cos({\pi s\over 2})\Ga(s)\ze(s)=\ze(1-s)\pi^s 2^{s-1}$
and the reality of $\ze$ on $\R,$
$$
\eqalignno{
{(1-z)\ze(z-1/2)\over  \ze'(z+1/2)} &=
{(1-z)(z-1/2)\Ga(z-1/2) \ze(z-1/2)\over \Ga(z+1/2) \ze'(z+1/2)}=\cr
{(1-z)(1/2-z)\over 2\pi}&{(1-\sin(\pi z))\over \cos(\pi z)}
{\ze(3/2-z)\over  \ze'(1/2-z)}\approx\cr
{y^2+1/2\over 2\pi}&(-i)\Bigl({\ze(3/2+z)\over  -\ze'(1/2+z)}\Bigr)^*,
&\eqnu
\label\etafun\eqnum* 
}
$$
where $z=iy,\ y\gg 0.$ Thus we need to check that 
$$
\eqalignno{
\rho(y)\equal&\Re(\ze({3\over 2}+iy)(\ze'({1\over 2}+iy))^*)=
\Re(i\ze({3\over 2}+iy)(\ze_y({1\over 2}+iy))^*)=\cr
&\Re(\ze({3\over 2}+iy)\Im(\ze_y({1\over 2}+iy))-
\Im(\ze({3\over 2}+iy)\Re(\ze_y({1\over 2}+iy))
&\eqnu
\label\rhore\eqnum* 
}
$$
is mainly positive for $\ze_y(1/2+iy)=\partial\ze(1/2+iy)/\partial y.$
The first term of the last difference dominates
(we will skip the explanation). Let us show that
it has a very good reson to be positive on the zeros $z=iy.$

In the first place, $\Re(\ze(3/2+iy))$ is always positive
([E], 6.6, pg.129). Then,
following  [E], 6.5,
$$
\eqalignno{
&\ze({1\over 2}+iy)=e^{-i\th(y)}Z(y),\ \th(y)=
\Im\log(\Ga({1\over 4}+i{y\over 2}))-{y\over 2}\log(\pi)
&\eqnu
\label\tzasnum\eqnum* 
}
$$
for the classical real-valued  function $Z(y)$ of $y\in \R_+.$
For instance, the zeros of  $\Im\ze(1/2+iy)=-Z(y)\sin(\th(y))$ 
are either  the zeros of $\ze(k+1/2)$ (coming from  $Z$) or 
the zeros of $\sin(\th(y)),$ the so-called Gram points. For small
$y,$ they appear alternately (pg.125, Gram's law), i.e. indeed
$\Im(\ze_y(1/2+iy))>0$
on the imaginary (+0.5) zeros of the zeta. Even when  Gram's law
fails the rule due to Rosser et al. states that excluding the zeros
from the so-called Gram blocks the sign is plus (8.4, pg.180). It holds
for at least  13.4 millions of zeros of $\ze$, many more than we need here.

However the above reduction and the last
argument are not 
sufficient to explain why the first ``negative'' zero for
the initial $\eta_+(z|2)$ (see (\ref\etaim)) is $z=1977.2714i.$
The sign of  $\Im(\ze_y(1/2+iy))$ becomes minus for the first time
at  $y=282.4651$. The next one is $y= 295.5733.$ Not very impressive.

There is something remarkable about $z=1977.2714i.$ Exactly at
this point, the computer program looking for the zeros of the zeta
has to become
more sophisticated. The reason is that this zero
is very close to the previous one
$z= 1977.1739i.$ If there is a  correlation between the zeros
with the linear $q$-approximations to the left and such pairs of zeros
then it could be of a certain
importance. The existence of nearly coincident zeros with  low extrema
of $Z$ between them ``almost''
contradicts  the Riemann hypothesis
([E], 8.4, pg.178) and, indeed,
``must give pause to even the most convinced believer''.

%
%
%
\vfil
\section{Symmetrization}
In conclusion, let us try to establish connections with 
Riemann's estimate of the number
of zeros in the critical strip
and the Riemann hypothesis. We will examine
the (anti)symmetrizations of the $q$-zeta functions
considered above with respect to $k\leftrightarrow  -k$.
It is the simplest way to ensure the functional equation.
Let us begin with $\ZZ^{im}_{+q}(k|d)$ from (\ref\zzimd)
and its meromorphic  continuation  $\ZZ^{neg}_{+q}(k|d)$   
to negative  $\Re k$ from (\ref\zzveed).
Since they behave differently for positive and negative $\Re k$,
one may hope to avoid unwanted zeros upon the symmetrization.

Let us  renormalize  $\ZZ^{im}_{+q}$ (cf. (\ref\zzbeta))
as follows:
$$
\eqalignno{
&\ga(k)\equal
\prod_{j=0}^\infty {(1-q^{j+2k})(1-q^{j+1})
\over (1-q^{j+k+1/2})(1-q^{j+k+1})},\cr
&\widehat{\ZZ}^{im}_{+q}(k|d)\equal
 \ga(k)\ZZ^{im}_{+q}(k|d),\ \widehat{\ZZ}^{neg}_{+q}(k|d)\equal
 \ga(k)\ZZ^{im}_{+q}(k|d),
&\eqnu\cr
\label\zzgamma\eqnum* 
&\lim_{a\to \infty}a^k\widehat{\ZZ}^{im}_{+q}(k|d) = \widehat{\ZZ}_+(k|d)
\equal k\sqrt{\pi}\ZZ_+(k|d),\ \Re k>0,\cr
&\lim_{a\to \infty}a^k\widehat{\ZZ}^{neg}_{+q}(k|d) = 
\widehat{\ZZ}_+(k|d),\ \Re k<0, \for \cr
&\ZZ_+(k|d)\ =\ d^{-1/2-k}(1-2^{1/2-k})\Ga(k+1/2)\ze(k+1/2).
}
$$
Cf. (\ref\zzstd), (\ref\zzbeta), and  (\ref\zzpstabl).
For $\Re k\ge 0$, we set 
$$
\eqalignno{
&\widehat{\ZZ}^{sym}_{+q}(k|d)\equal 
 a^{k}\widehat{\ZZ}^{im}_{+q}(k|d) + 
 a^{-k}\widehat{\ZZ}^{neg}_{+q}(-k|d), &\eqnu \cr
\label\zzgsym\eqnum*
&\lim_{a\to \infty}\widehat{\ZZ}^{sym}_{+q}(k|d) = 
\widehat{\ZZ}^{sym}_+(k|d)
\equal k\sqrt{\pi}(\ZZ_+(k|d)- \ZZ_+(-k|d)).
}
$$

The function $\widehat{\ZZ}^{sym}_{+q}$
can be naturally extended to a 
symmetric (even) 
meromorphic function defined for all $k\in \C$: 
$\widehat{\ZZ}^{sym}_{+q}(-k|d) =$ $\widehat{\ZZ}^{sym}_{+q}(k|d)$.
It has zeros at $k\in 2\pi ia\Z$. Its poles modulo $2\pi a i \Z$ 
form (see (\ref\sbul)) the set
$$
\eqalignno{
&\widehat{\p}=(-\La)\cup\La\cup \{1/2+\Z\},\ 
\La= \{-2\xi_l -\Z_+\}\cup
 \{-2\xi_l^* -\Z_+ \}.
&\eqnu
\label\sbuls\eqnum*
}
$$
They are all simple for generic $a$.

Given $0<\ka<1$ such that $\{\ka+i\R+\Z\}\cap\widehat{\p}=\emptyset$,
$$
\eqalignno{
& \lim_{r\to \infty} a^{r-\ka}
\widehat{\ZZ}^{sym}_{+q}(\ka-r+i\nu)\ =\ a^{i\nu} \psi_+,\cr
&\where\psi_+\ =\ \lim_{\Re k\to\infty}\ZZ^{im}_{+q}(k)>0, \
r\in \N.
&\eqnu
\label\psisym\eqnum* 
}
$$
Here we employ Theorem \ref\ZZKLIM and 
(\ref\zzlimplus) and use that  $(\ga\be^{-1})(\ka-r+\nu i)$ 
(see (\ref\zzbeta)) tends to zero as $\N\ni r\to-\infty$.  

Because of the appearance of the $a$-factors the
function $\widehat{\ZZ}^{sym}_{+q}(k)$ is
not $2\pi a i$-periodic. However we may assume that
$a^{2\pi ai}=1$, i.e. 
$$ai \log(a)\ =\ M \for M\in\Z.$$
Then it becomes $2Ti$-periodic for $T\equal \pi a$ and
we can readily calculate the 
number of its zeros minus the number of poles modulo $2Ti$
in the limit of big $r$. Moreover, because of  the reality
of $\widehat{\ZZ}^{sym}_{+q}(k)$ on the real axis, if
 $z$ is a zero (a pole) so is $2Ti-z.$
Hence we can find this  number 
for the half-period too.
  
\proclaim {Corollary}
Let $\ka$ be as above,  $\widehat{O}_r$ the number of zeros of
$\widehat{\ZZ}^{sym}_{+q}(k)$  in the strip  
$$\s_r\equal \{\ka-r\le \Re k \le r-\ka,\  0 \le \Im k  \le T\},$$
$\widehat{P}_r$ the number of poles. We calculate 
zeros(poles) with multiplicities 
and with the coefficient $1/2$ if they belong to the boundary.
Then 
$$
\eqalignno{
& \widehat{N}_\infty\equal  \widehat{O}_r-\widehat{P}_r \ =\ 
 {T\over \pi}(\log{T\over\pi})
\for \N\ni r \gg 0.
&\eqnu
\label\nopr\eqnum* 
}
$$
\label\NOPR\theoremnum*
\endproclaim

Formula (\ref\nopr) can be considered as a certain counterpart
of the celebrated formula (see [E, 6.7])
$$
\eqalignno{
& N(T)\ = \ {T\over 2\pi}(\log{T\over 2\pi}-1)+O(\log T)
&\eqnu 
\label\noclz\eqnum* 
} 
$$
for the number of the zeros $z$ 
of the zeta with $0<\Im z<T$  (in the critical strip).
The total  number of imaginary  zeros of 
$\widehat{\ZZ}_{+}^{sym}(k|d)$ such that  $0<\Im z<T$ results
from this formula modulo the Riemann hypothesis:
$$
\eqalignno{
&\widehat{N}(T|d)={T\over\pi}(\log T-1-\log(2d))+O(\log T). 
&\eqnu
\label\nopzes\eqnum* 
}
$$
All zeros of $\ze(k+1/2)$ are those of
 $\widehat{\ZZ}_{+}^{sym}.$ However the latter may have extra zeros.
According to our analysis their number is finite for  
$d\le 1$ and is not bigger than $(T/\pi)\log 2$  for $d>1$
For big $\Im k$
they may appear only near the zeros of $(1-2^{1/2\pm k}).$
Anyway the non-imaginary zeros do not change
(\ref\nopzes) too much.

Concerning $\widehat{\ZZ}_{+q}^{sym}(k|d)$,
first, it must have exactly one  zero 
approaching  each pole from the set $\{\pm k=1/2+2\Z_+\}$ as $q\to 1$.
Indeed,
the limit  $\widehat{\ZZ}_{+}^{sym}(k|d)$ has neither poles
nor zeros  there.
It is likely that the same holds for the remaining poles at
$\{\pm k=-1/2+2\Z_+\}.$ At least it is true
for sufficiently big $|\Re k|$ because of the stable 
 $\Z$-periodicity of $\ZZ_{+}^{neg}(k|d).$ 
Numerical experiments confirm this: the ``symmetric''
$q$-zeros practically
coincide with those  from  (\ref\stzero) as $a=1=d$. 

Second, there is a tendency for the poles from $\La$ to be moving
accompanied by the zeros.  
We cannot prove that this always happens.
Anyway the number of $q$-poles equals the number
of $q$-zeros inside a closed curve in the form
$\c_a=\sqrt{a}\c_1$ for a given $\c_1$
if  
$$|a^{k}\widehat{\ZZ}^{im}_{+q}(k|d)| > 
|a^{-k}\widehat{\ZZ}^{neg}_{+q}(k|d)|\for k\in \c_a,
$$
assuming that $\Re \c_1>0$ and $\widehat{\ZZ}^{im}_{+q}(k|d)$
has no zeros inside $\c_a.$
Generally speaking, it holds.
However it is not clear at the moment whether
this can be used for  quantitative
estimates. 

Anyway if it is true (or almost true)
that the poles of  $\widehat{\ZZ}_{+q}^{sym}(k|d)$ from $\La$
are compensated
by the zeros in the limit $a\to \infty,$ then we come to the following
conjecture. 

\proclaim {Conjecture} Assuming that $d\ge \pi/(2e),$
the $q$-deformation $z_q$ of  any imaginary zero $z=iy (y>0)$ of 
$\widehat{\ZZ}^{sym}_{+}(k|d),$ namely a zero
of $\widehat{\ZZ}^{sym}_{+q}(k|d)$ approaching $z$
in the limit $q\to 1$, exists and  remains in the horizontal strip
$0< \Im z_q < y + C$ for all $0<q<1$ and a proper constant $C>0$.
The zeros of $\widehat{\ZZ}^{sym}_{+q}(k|d)$ in the box
$$
\eqalignno{
&K_a(c)\ =\ \{-2\sqrt{\pi a\over 2d}+c < \Re k<2 \sqrt{\pi a\over 2d}-c,
\ 0<\Im k<\pi a\}
&\eqnu
\label\constrip\eqnum* 
}
$$
are $q$-deformations of the classical ones
for a constant $c>0$.
\label\CONN\theoremnum*
\endproclaim

The first part of the conjecture for $d= \pi/(2e)$ 
results in the  estimate  
$$N(T)\ < \ {T\over 2\pi}(\log{T\over 2\pi}-1)+O(\log T),$$
speculating that the number of poles of $\ZZ^{sym}_{+q}(k|d)$ for big $a$
with $0<|\Re k|<t$ which have more than one neighboring
zero is $O(\log t)$ at most. As $a\to 0$ (i.e. $q=\exp(-1/a)\to 0$),
the $q$-deformations must approach the real line. So eventually
they go down for small $a$. The box $K_a(0)$ in the second
part does not contain the poles of $\ZZ^{sym}_{+q}.$  Thus we
expect that once a $q$-zero gets distant enough from the $q$-poles it
always approaches a classical zero in the limit.

\vskip 5pt   {\bf Variants.}
Now let us modify the imaginary $q$-zeta to avoid the gamma-factors. 
For $\Re k\ge 0$, we set 
$$
\eqalignno{
&\hze_{+q}(k|d)\equal 
  4\si(k)\ZZ^{im}_{+q}(k|d) - 
 4\si(-k)\ZZ^{neg}_{+q}(-k|d),\cr
&\where \si(k)\equal
\prod_{j=0}^\infty (1-q^{k+1/2+j})/( 1-q^{j+1}),
&\eqnu \cr
\label\zzsym\eqnum*
&\lim_{a\to \infty}\hze_{+q}(k|d)\ =\ \hze_+(k|d)
\equal \ze_+(k|d)- \ze_+(-k|d),\cr
&\ze_+(k|d)\ =\ (d/4)^{-1/2-k}(1-2^{1/2-k}) \ze(k+1/2).
}
$$
The function $\hze_{+q}$ can be extended to an odd
meromorphic $2\pi a i$-periodic function with the poles 
at $\{\pi ai +\Z/2\}.$ Since $(\si\be^{-1})(k-r+\nu i)\to 0$
as $ r\to 0$ (see (\ref\psisym)), 
$$
\eqalignno{
& \lim_{r\to \infty} \hze_{+q}(k|d)(\ka-r+i\nu)\ =\ 4 \psi_+,\ r\in \N.
&\eqnu
\label\psiasym\eqnum* 
}
$$
Hence the (stable)
numbers of  zeros of  $\hze_{+q}(k|d)$ and poles
on the half-period
$0\le \Im k< T=\pi a$ coincide. So
there is no relation to the Riemann
estimate in the considered case, 
which can be explained as follows.
First, the classical function   $\hze_{+}(k|d)$
has infinitely many real zeros approaching the
points from $\pm\{1/2+2\Z_+\}$ for big $|\Re k|.$
It can also have zeros which are neither real nor imaginary 
(but finitely many). Second, the
poles of $\hze_{+}(k|d)$ on the line $\Im k=\pi a$ may ``loose''
the neighboring zeros as $a\to \infty$. 
In line with Conjecture \ref\CONN, such lost $q$-zeros could be
a significant  source of the classical ones.
Hopefully it may be seen numerically. 

\vskip 5pt
It is  instructional to get a $q$-deformation
of the modified zeta
$$\tze(k)=\pi^{-(k/2+1/4)}\Ga ({k\over 2}+{1\over 4}) \ze(k+{1\over 2})$$
satisfying the plane functional equation $\tze(k)=\tze(-k).$ 
As always, we need to switch to its plus-variant to ensure the convergence
for $\Re k\le 1/2$:
$$
\eqalignno{
& \tze_{+}(k|d)\ =\ p(k|d)\tze(k),\  p(k|d)\ =\cr
&d^{-1/2-k}(1-2^{1/2-k})-d^{-1/2+k}(1-2^{1/2+k}).
&\eqnu
\label\zesym\eqnum* 
}
$$
All zeros of the  function $ p(k|d)$ are imaginary for $d\ge 1/2.$ 
The number of imaginary zeros
in the interval  $0<\Im k<T$ equals $|(T/\pi)\log(2d)|$
up to $-1\le \ep\le 1$ for any $d>0.$
Let us assume that $d>1/2.$

We set
$$
\eqalignno{
&\al(k)\equal \ga(k)
\prod_{j=0}^\infty { 1-q^{j+1}\over 1-q^{k/2+1/4+j} },\ \
(4\pi)\,\tze_{+q}(k|d)\equal
&\eqnu
\label\zzals\eqnum* \cr
& ({a\over 16\pi})^{k/2-3/4}\al(k)\ZZ^{im}_{+q}(k|d) -
({a\over 16\pi})^{-k/2-3/4}\al(-k)\ZZ^{neg}_{+q}(-k|d).
}
$$
Then $\lim_{a\to \infty}\tze_{+q}(k|d)\ =\ \tze_+(k|d).$

The function $(\al\be^{-1})(\ka-r+\nu i))$ tends to zero 
as $r\to \infty$ too, so we can calculate the stable (see above) 
difference $ \widetilde{N}_\infty$ between 
the total number of zeros and poles
on the ``half-period'' for $\tze_{+q}(k|d)$. 
We may assume that
$(a/(16\pi))^{\pi ai/2}=1$, i.e. 
$$(\pi ai) \log(a/(16\pi)) =M \for M\in\Z,$$
to provide the $2\pi ai$-periodicity of  $\tze_{+q}.$
Then $ \widetilde{N}_\infty\ =\ M.$
Let us express it in terms of $T=\pi a:$
$$
\eqalignno{
& \widetilde{N}_\infty(T)\ =\ M\ = \   
{T\over 2\pi}(\log{T\over 16\pi}-\log\pi).
&\eqnu
\label\noprs\eqnum* 
}
$$
Here the leading term coincides  with that from
(\ref\noclz) again. More exactly, 
the  number of zeros of 
$\tze_{+}(k|d)$ such that  $0<\Im z<T$ equals
$$
\eqalignno{
& \widetilde{N}(T|d)={T\over 2\pi}(\log{T\over 2\pi}-1+2\log(2d))+O(\log T). 
&\eqnu
\label\nopzess\eqnum* 
}
$$

The zeros
of $\tze_{+q}$ either a) approach its  real poles at
$\{\pm k=1/2+2\Z_+\}$, or 
b) go together with the corresponding poles from $\La$
(see (\ref\sbuls)) upwards, or
c) tend to the  zeros of zeta in the critical strip.
We may conjecture that 
a $q$-deformation $z_q$ of  a zero $z=iy (y>0)$ of $\ze(k+1/2)$ 
remains in the horizontal strip
$0< \Im z_q < y + C_d y/\log y$ for all $0<q<1$
and a constant $C_d>0$ depending on $d$.
One can also speculate that
$C_d\sim \log(32\pi d^2/e)$ comparing  
the formulas for $\widetilde{N}(T|d)$ and $\widetilde{N}_\infty(T)$
in more detail.

\vskip 5pt
The above estimates can be extended to the
case of  the standard zeta. The factors $(1-2^{1/2-k})$
will disappear from the formulas. The calculation 
of the stable number of zeros minus poles remains the same. However the 
minus-counterparts $\widehat{\ZZ}^{sym}_q,
\ \hze_{q},$ and $\tze_{q}$
of the functions considered above
diverge in the critical
strip. So the relations to the classical zeros have to be
discussed under the constraint
$\log(a)=o(|\Ga(k)^2\tan(\pi k)|^{-1/2}).$
The conjecture and  other considerations can be also
extended to the $q$-deformations 
of the Dirichlet $L$-function (see (\ref\zzdiri)),
where the convergence holds for all $k$ except for the poles.

\vskip 5pt  
{\bf Towards Riemann $q$-hypothesis.}
The conjecture states nothing about the pure
imaginary zeros of $\widehat{\ZZ}^{sym}_{+q}(k|d).$ 
It is more convenient to discuss this problem for the
sharp zeta-function. For the sake of concreteness, we
assume that $d=1$ and will $q$-deform the $\hze$
from (\ref\zzsym) only. A similar behavior is  expected for 
the sharp-variants of the remaining two symmetric $q$-zeta
functions considered above, their minus-counterparts
(for $a$ which are not too big), and the sharp
Dirichlet $L$-functions.

We set
$$
\eqalignno{
&\ \ \ \hze^\sharp_{+q}(k)\equal 
  \chi(k)^{-1}\ZZ^{\sharp}_{+q}(k) - 
 \chi(-k)^{-1}\ZZ^{\sharp}_{+q}(-k),&\eqnu \cr
\label\hzesha\eqnum*
&\chi(k)\equal
\GG_q^{\sharp}(k) = -{\pi a\over 2}q^{k^2}(1-q^{-k})
\prod_{j=0}^\infty
{(1-q^{j+k}) (1-q^{j/2+1/2})\over
 (1-q^{j+2k})(1-q^{j+1})^2 }\times\cr
& \prod_{j=1}^\infty
(1-q^{j-k})(1-q^{j+k})(1+q^{j/2-1/4+k/2})(1+q^{j/2-1/4-k/2}),\cr
&\lim_{a\to \infty}\hze_{+q}(k)\ =\ \hze_+(k)\ =
\equal \ze_+(k)- \ze_+(-k).
}
$$
Thus we just divided the sharp-integral for the
kernel $(1+q^{-x^2})$ by that for the Gaussian $q^{x^2}.$
This very choice of the normalization is not
too important. However in this form the definition can be naturally
extended to the multi-dimensional case. 
Thanks to Theorem \ref\GG and Theorem \ref\ZZSP, the function
$\hze^\sharp_{+q}(k)$ 
is regular in the horizontal strip 
$$K_a^{\sharp}(c)\equal \{-\sqrt{2\pi a}+c <
\Im k< \sqrt{2\pi a}-c\}$$
for  $c\ge  0$.
 
\proclaim {Conjecture} 
Let us fix $q=\exp(-1/a)$ for $a>0.$ Then
all zeros  of 
$\hze_+(k)$ inside $K_a^{\sharp}(0)$ have unique
$q$-deformations, i.e. 
the zeros of $\hze^\sharp_+(k)$ convergent to the corresponding
classical ones in the limit, which belong to  $K_+^{\sharp}(-c)$ for
a constant $c>0.$ Vice versa, all zeros
of $\hze^\sharp_+(k)$ inside $K_a^{\sharp}(c)$ 
 are such deformations for proper $c>0.$
The $q$-deformations of the imaginary zeros of
$\hze_+(k)$ in $K_a^{\sharp}(0)$ are imaginary.
\label\CONR\theoremnum*
\endproclaim

The  sharp-counterpart of the Riemann hypothisis is the claim 
that all zeros $z_q$ of $\hze^\sharp_{+q}(k)$
in the critical strip $\{|\Re k|<1/2\}$ such that
$0<\Im z_q< \sqrt{2\pi a}-c$  for a certain
constant $c>0$ are imaginary.
Confirmations of the conjecture are entirely numerical.
No theoretical approach is known at the moment.

\comment
Shi[k_Complex]:=
Kzi[k+1.]-(1-Exp[-(2k+1)/a])*(1+Exp[-k/a])*(1-Exp[-(k+1)/a])^(-1)*Kzi[k];
Shh[k_Complex]:=Exp[-(1.5+k)/a]*((1-Exp[-(2k+2)/a])*
(1-Exp[-(2k+3)/a]))^(-1)*Kzis[k+2.];
THM: ALWAYS when ss(iR+1) and ss(iR-1) do not have poles between,
v(0+1)=r^4 (v+v^(-1) is 2 at 0 and increases till r^4+r^(-4) at 1,-1).
with respect to z.
SHI+SHH==0!!!
If E(x) has no poles from -u<Re[x]<u then Kzi[k] is analytic
for Rek>-2u with poles only at -.5, -1.5 etc.
\endcomment

%
%
%
%
%
%
\AuthorRefNames [BG]
\references
\vfil

[A] edited by \name{G. E.  Andrews} [et al],
{ Ramanujan Revisited}, Academic Press Inc.
(1988).

[AI]
\name{R. Askey}, and \name{M.E.H. Ismail},
{ A generalization of ultraspherical polynomials},
in Studies in Pure Mathematics (ed. P. Erd\"os),
Birkh\"auser (1983), 55--78.

[AW]
\bibline, and \name{J. Wilson},
{ Some basic hypergeometric orthogonal polynomials
that generalize Jacobi polynomials},
Memoirs AMS {319} (1985).

[BI] 
\name{B.M. Brown}, and \name{M.E.H. Ismail},
{ A right inverse of the Askey-Wilson operator},
Proceedings of AMS {123} (1995), 2071--2079.

[C1]
\name{I. Cherednik},
{ Difference Macdonald-Mehta conjecture},
IMRN {10} (1997), 449--467.

[C2]
\bibline,
{ Inverse Harish-Chandra transform and difference ope\-ra\-tors},
IMRN {15} (1997), 733--750.

[C3]
\bibline, 
{ Macdonald's evaluation conjectures and
difference Fourier transform},
Inventiones Math. {122} (1995),119--145.

[E]
\name{H.M. Edwards},
{ Riemann's zeta function},
Academic Press, New York and London (1974). 

[H]
\name{G.J. Heckman},
{ An elementary approach to the hypergeometric shift operators of
Opdam}, Invent.Math. {103} (1991), 341--350.

[HO]
\bibline, and \name{E.M. Opdam}, 
{ Root systems and hypergeometric functions I},
Comp. Math. {64} (1987), 329--352. 

[He]
\name {S. Helgason},
{ Groups and geometric analysis}, 
Academic Press, New York (1984). 

[J]
\name{M.F.E. de Jeu},
{ The Dunkl transform}, Invent. Math. {113} (1993), 147--162.

[M1]
\name{I.G. Macdonald}, 
{ Some conjectures for root systems},
SIAM J.Math. Anal. {13}:6 (1982), 988--1007.

[M2]
\name{I.G. Macdonald}, 
{ A new class of symmetric functions },
Publ.I.R.M.A., Strasbourg, Actes 20-e Seminaire Lotharingen,
(1988), 131--171 .

[Mo]
\name{D.S. Moak},
{ The $q$-analogue of Stirling formula},
Rocky Mountain J. Math. {14} (1984), 403--413.

[O]
\name{E.M. Opdam}, 
{ Some applications of hypergeometric shift
operators}, Invent. Math. {98} (1989), 1--18.

[R]
\name {B. Riemann},
{ Gesammelte Werke}, Teubner, Leipzig (1892),
Dover Books, New York (1953).

[Sa]
\name {J. Satoh},
{ $q$-Analogue of Riemann's $\zeta$-function and $q$-Euler numbers},
J. Number Theory {31} (1989), 346--362.

[UN]
\name{K. Ueno}, \and {M. Nishizawa},
{ Quantum groups and zeta-functions}, 
 Proceedings of the 30-th Karpatz Winter School
 ``Quantum Groups: Formalism and Applications'' (1995), 115--126
 (Polish Scientific Publishers PWN).

\endreferences

\bye